 \documentclass[12pt]{article}

 \usepackage{mathrsfs}
 \usepackage{amssymb}

 \def\proclaim{\medskip\noindent}\def\endproclaim{\medskip}

 \def\beqlb{\begin{eqnarray}}\def\eeqlb{\end{eqnarray}}
 \def\beqnn{\begin{eqnarray*}}\def\eeqnn{\end{eqnarray*}}

 \def\qed{\hfill$\Box$\medskip}

 \def\<{\langle}\def\>{\rangle}

 \def\n{~}

 \def\d{\mbox{\rm d}}\def\e{\mbox{\rm e}}

 \def\IR{I\!\! R}

 \def\bE{\mbox{\boldmath $E$}}\def\bP{\mbox{\boldmath $P$}}
 \def\bQ{\mbox{\boldmath $Q$}}\def\bU{\mbox{\boldmath $U$}}

\begin{document}

\

\noindent{Published in: {\it Theory of Probability and its
Applications} {\bf 46} (2002), 274--296.}

\bigskip\bigskip
\bigskip

\centerline{\Large\bf SKEW CONVOLUTION SEMIGROUPS}

\bigskip

\centerline{\Large\bf AND RELATED}

\bigskip

\centerline{\Large\bf IMMIGRATION PROCESSES\,\footnote{Supported
by the National Natural Science Foundation of China
(No.10131040).}}

\bigskip

\centerline{Zeng-Hu LI}

\medskip

\centerline{\it Department of Mathematics, Beijing Normal University,}

\centerline{\it Beijing 100875, People's Republic of China}

\centerline{\it E-mail: lizh@email.bnu.edu.cn}

\bigskip\bigskip

{\narrower{\narrower

\noindent{\bf Abstract.} A special type of immigration associated
with measure-valued branching processes is formulated by using
skew convolution semigroups. We give characterization for a
general inhomogeneous skew convolution semigroup in terms of
probability entrance laws. The related immigration process is
constructed by summing up measure-valued paths in the Kuznetsov
process determined by an entrance rule. The behavior of the
Kuznetsov process is then studied, which provides insights into
trajectory structures of the immigration process. Some well-known
results on excessive measures are formulated in terms of
stationary immigration processes.

\bigskip


\bigskip


\par}\par}

\bigskip\bigskip

\section{Introduction}

\setcounter{equation}{0}

Let $E$ be a Lusin topological space, i.e., a homeomorph of a
Borel subset of a compact metric space, with the Borel
$\sigma$-algebra ${\cal B}(E)$. Let $B(E)$ denote the set of
bounded ${\cal B}(E)$-measurable functions on $E$, and $B(E)^+$
the subspace of $B(E)$ of non-negative functions. We denote by
$M(E)$ the space of finite measures on $(E,{\cal B}(E))$ endowed
with the topology of weak convergence. For $f\in B(E)$ and $\mu\in
M(E)$, write $\mu(f)$ for $\int_E f\d\mu$. Suppose that
$(P_t)_{t\ge0}$ is the transition semigroup of a Borel right
process $\xi$ with state space $E$ and $\phi(\cdot,\cdot)$ is a
branching mechanism given by
 \beqlb
\phi(x,z) = b(x)z +
c(x)z^2 + \int_0^\infty(\e^{-zu}-1+zu) m(x,\d u),
\quad x\in E,\ z\ge0,
\label{1.1}
 \eeqlb
where $b\in B(E)$, $c\in B(E)^+$ and $(u\land u^2) m(x,\d u)$ is a
bounded kernel from $E$ to $(0,\infty)$. Then for each $f\in
B(E)^+$ the evolution equation
 \beqlb
V_tf(x) + \int_0^t\d s\int_E \phi(y,V_sf(y))
P_{t-s}(x,\d y) = P_tf(x),
\quad t\ge0,\ x\in E,
\label{1.2}
 \eeqlb
has a unique solution $V_tf\in B(E)^+$, and there is a Markov
semigroup $(Q_t)_{t\ge0}$ on $M(E)$ such that
 \beqlb
\int_{M(E)} \e^{-\nu(f)} Q_t(\mu,\d\nu) =
\exp\left\{-\mu(V_tf)\right\} \label{1.3}
 \eeqlb
for all $t\ge0$, $\mu\in M(E)$ and $f\in B(E)^+$. A Markov process
$X$ having transition semigroup $(Q_t)_{t\ge0}$ is called a {\it
Dawson-Watanabe superprocess} with parameters $(\xi,\phi)$. Under
our hypotheses, $X$ has a Borel right realization; see Fitzsimmons
\cite{F88} and \cite{F92}. The $(\xi,\phi)$-superprocess is a
mathematical model for the evolution of a population in some
region; see e.g. Dawson \cite{D92} and \cite{D93}. If we consider
a situation where there are additional sources of population from
which immigration into the region occurs during the evolution, we
need to introduce branching processes with immigration. This type
of modification is familiar from the branching process literature;
see e.g. \cite{AN72}, \cite{DI78}, \cite{D91}, \cite{KW71},
\cite{L92b} and \cite{S90}.

A class of measure-valued immigration processes were formulated in
Li \cite{L96a} as follows. Let $(N_t)_{t\ge0}$ be a family of
probability measures on $M(E)$. We call $(N_t)_{t\ge0}$ a {\it
skew convolution semigroup} associated with $X$ or $(Q_t)_{t\ge0}$
provided
 \beqlb
N_{r+t} = (N_rQ_t)*N_t,
\quad r,t\ge0,
\label{1.4}
 \eeqlb
where ``$*$'' denotes the convolution operation. The relation
(\ref{1.4}) is satisfied if and only if
 \beqlb
Q^N_t(\mu,\cdot):= Q_t(\mu,\cdot)*N_t, \quad t\ge0,\ \mu\in M(E),
\label{1.5}
 \eeqlb
defines a Markov semigroup $(Q^N_t)_{t\ge0}$ on $M(E)$. A Markov
process is called an {\it immigration process} associated with $X$
if it has transition semigroup $(Q^N_t)_{t\ge0}$. The intuitive
meaning of the immigration process is clear from (\ref{1.5}), that
is, $Q_t(\mu,\cdot)$ is the distribution of descendants of the
people distributed as $\mu\in M(E)$ at time zero and $N_t$ is the
distribution of descendants of the people immigrating to $E$
during the time interval $(0,t]$. Clearly, (\ref{1.5}) gives the
general formulation for the immigration independent of the inner
population.

Needless to say, most of the theory of Dawson-Watanabe
superprocesses carries over to their associated immigration
processes and could be developed by techniques very close to those
in \cite{D92} and \cite{D93}. It is interesting, however, that the
immigration processes have many additional structures, as might be
expected from (\ref{1.4}) and (\ref{1.5}). Note that (\ref{1.5})
is quite similar to the construction of L\'evy's transition
semigroup from a {\it usual} convolution semigroup. It is
well-known that a convolution semigroup on the Euclidean space is
uniquely determined by an infinitely divisible probability
measure. As shown in Li \cite{L96a}, the skew convolution
semigroup may be characterized in terms of an infinitely divisible
probability entrance law. Therefore, the immigration process may
be regarded as a generalized form of the celebrated L\'evy
process. Other examples of immigration processes are squares of
Bessel diffusions and radial parts of Ornstein-Uhlenbeck
diffusions; see \cite{KW71} and \cite{SW73}. The above formulation
also includes new kinds of processes; see \cite{L96b}, \cite{L98a}
and \cite{L98b}.

In the immigration models studied before, authors usually assumed
that the immigrants came to $E$ according to a random measure on
$\IR\times E$. The scenes are not quite clear under the skew
convolution semigroup formulation. This weak point has in fact
motivated the present work. The main purpose of this paper is to
give interpretations of the skew convolution semigroups by
constructing and analyzing the trajectories of the related
measure-valued immigration processes. We first prove that a
general inhomogeneous skew convolution semigroup may be decomposed
into three components which involve respectively a countable
family of entrance laws, a countable family of closed entrance
laws, and a continuum family of infinitely divisible probability
entrance laws together with a diffuse measure on the index space.
Then we give a construction for the immigration process defined
above by picking up measure-valued paths with random times of
birth and death. Our construction is based on the observation that
any skew convolution semigroup is determined by a continuous
increasing measure-valued path $(\gamma_t)_{t\ge0}$ and an
entrance rule $(G_t)_{t\ge 0}$. This fact yields a natural
decomposition of the immigration into two parts; the deterministic
part represented by $(\gamma_t)_{t\ge0}$ and the random part
determined by $(G_t)_{t\ge0}$. The latter is usually an
inhomogeneous immigration process and can be constructed by
summing up paths $\{w_t: \alpha<t<\beta\}$ in the associated
Kuznetsov process; see Kuznetsov \cite{K74}. By analyzing the
asymptotic behavior of the paths $\{w_t: \alpha<t<\beta\}$ near
the birth time $\alpha= \alpha(w)$, we show that almost all these
paths start propagation in an extension $E^T_D$ of the underlying
space, including those growing up at points in this space from the
null measure. Those combined with our construction of the
immigration process give a full description of the phenomenon. In
some special cases, the infinitely divisible entrance law for the
$(\xi,\phi)$-superprocess corresponds to a $\sigma$-finite
entrance law and the associated immigration process can be
constructed using a homogeneous path-valued Poisson random process
whose characteristic measure is the Markov measure determined by
the entrance law. The construction has been proved useful in
studying the immigration processes; see e.g. \cite{L96b},
\cite{LS95} and \cite{S90}. As additional applications of the
construction, we give formulations of some well-known results on
excessive measures in terms of stationary immigration processes.

The paper is organized as follows. The next section contains some
preliminaries. In section~3, we prove the decomposition theorem
for inhomogeneous skew convolution semigroups. The construction of
immigration processes using Kuznetsov processes is given in
section~4. Almost sure behavior of the Kuznetsov processes is
studied in section~5. In section~6, we discuss stationary
immigration processes determined by excessive measures.

\section{Preliminaries}
\setcounter{equation}{0} Recall that $M(E)$ is the space of finite
Borel measures on the Lusin topological space $E$. It is
well-known that $M(E)$ endowed with the weak convergence topology
is also a Lusin space. Let $M(E)^\circ=M(E)\setminus\{0\}$, where
$0$ denotes the null measure. For a probability measure $F$ on
$M(E)$ we define its Laplace functional by
 \beqlb
L_F(f) := \int_{M(E)} \e^{-\nu(f)} F(\d\nu),
\quad f\in B(E)^+,
\label{2.1}
 \eeqlb
which determines $F$ uniquely. It is well-known that $F$ is
infinitely divisible if and only if
 \beqlb
L_F(f) = \exp\bigg\{-\eta(f) -
\int_{M(E)^\circ}\left(1-\e^{-\nu(f)} \right)H(\d\nu)\bigg\},
\quad f\in B(E)^+, \label{2.2}
 \eeqlb
where $\eta\in M(E)$ and $[1\land\nu(E)] H(\d\nu)$ is a finite
measure on $M(E)^\circ$. See e.g. Kallenberg \cite{K75}. We write
$F=I(\eta,H)$ if $F$ is determined by (\ref{2.2}).

Suppose that $X = (W, {\cal G}, {\cal G}_{r,t}, X_t, \bQ_{r,\mu})$
is a Markov process in $M(E)$ with transition semigroup
$(Q_{r,t})_{r\le t}$. Let
 \beqlb
V_{r,t}f(x) :=
-\log\int_{M(E)}\e^{-\nu(f)}Q_{r,t}(\delta_x,\d\nu), \quad r\le
t,\ x\in E,\ f\in B(E)^+, \label{2.3}
 \eeqlb
where $\delta_x$ denote the unit mass concentrated at $x\in E$. In
this paper, we always assume that, for each $r\le t$ and $f\in
B(E)^+$, the function $V_{s,t}f(x)$ of $(s,x)$ restricted to
$[r,t]\times E$ is bounded and measurable. We also assume that
$V_{r,t}f(x)$ is right continuous in $t\ge r$ for $f\in C(E)^+$,
continuous functions in $B(E)^+$. The process $X$ is called a {\it
measure-valued branching process} (MB-process) if its transition
semigroup satisfies
 \beqlb
\int_{M(E)} \e^{-\nu(f)} Q_{r,t}(\mu,\d\nu) =
\exp\left\{-\mu(V_{r,t}f)\right\}, \quad r\le t,\ f\in B(E)^+.
\label{2.4}
 \eeqlb
Under this hypothesis, $Q_{r,t}(\mu,\cdot)$ is infinitely
divisible and $(V_{r,t}) _{r\le t}$ form a family of operators on
$B(E)^+$ satisfying $V_{r,s}V_{s,t} =V_{r,t}$ for all $r\le s\le
t$, which is called the {\it cumulant semigroup} of $X$. See e.g.
Silverstein \cite{S69} and Watanabe \cite{W68}. The
$(\xi,\phi)$-superprocess defined in the introduction is a special
form of the MB-process.

Suppose that $(Q_{r,t})_{r\le t}$ is the transition semigroup of
an MB-process and $(N_{r,t}) _{r\le t}$ is a family of probability
measures on $M(E)$. We call $(N_{r,t}) _{r\le t}$ a {\it skew
convolution semigroup} (SC-semigroup) associated with
$(Q_{r,t})_{r\le t}$ if
 \beqlb
N_{r,t} = (N_{r,s}Q_{s,t})*N_{s,t},
\quad r\le s\le t.
\label{2.5}
 \eeqlb

\proclaim{\bf Theorem\n{2.1}.} {\sl The equation (\ref{2.5}) is
fulfilled if and only if
 \beqlb
Q^N_{r,t}(\mu,\cdot):= Q_{r,t}(\mu,\cdot)*N_{r,t}, \quad r\le t,\
\mu\in M(E), \label{2.6}
 \eeqlb
defines a Markov semigroup $(Q^N_{r,t})_{r\le t}$ on $M(E)$.}
\endproclaim

\noindent{\it Proof.} Let $(Q_{r,t})_{r\le t}$ and
$(Q^N_{r,t})_{r\le t}$ be given by (\ref{2.4}) and (\ref{2.6}),
respectively. It is not difficult to check that (\ref{2.5}) is
equivalent to the Chapman-Kolmogorov equation
 \beqnn
Q^N_{r,t}(\mu,\cdot) =
\int_{M(E)}Q^N_{r,s}(\mu,\d\nu)Q^N_{s,t}(\nu,\cdot), \quad r\le
s\le t,\ \mu\in M(E),
 \eeqnn
from which the assertion follows. \qed

Let $T$ be an interval and $(N_{r,t})_{r\le t}$ be an SC-semigroup
associated with $(Q_{r,t})_{r\le t}$. If $\{Y: t\in T\}$ is a
Markov process having transition semigroup $(Q^N_{r,t})_{r\le t}$
defined by (\ref{2.6}), we call it an {\it immigration process}
associated with $X$. (Of course, SC-semigroups and immigration
processes can also be formulated for some more general classes of
Markov processes with state spaces possessing semigroup
structures.)

It is known that a metric $\varrho$ can be introduced into $E$ so
that the Borel $\sigma$-algebra on $E$ induced by $\varrho$
coincides with its original Borel $\sigma$-algebra; see e.g. Cohn
\cite[p275]{C80}. We write $M(E_\varrho)$ for the set $M(E)$
furnished with the topology of weak convergence on $(E,\varrho)$.
Then $M(E_\varrho)$ is locally compact and metrizable. Let
$D(E_\varrho)^+$ be a countable dense subset of the space of {\it
strictly positive} continuous functions on $(E,\varrho)$.

\proclaim{\bf Lemma\n{2.2}.} {\sl Let $\{F_n: n=1,2,\cdots\}$ be a
sequence of probabilities on $M(E)$. If the limit
 \beqlb
L(f) := \lim_{n\to\infty} L_{F_n}(f), \quad f\in D(E_\varrho)^+,
\label{2.7}
 \eeqlb
exists and $L(f)\to 1$ as $f\to 0$, then there is a probability
measure $F$ on $M(E)$ such that $L_F(f) = L(f)$ for all $f\in
D(E_\varrho)^+$. Moreover, if each $F_n$ is infinitely divisible,
so is $F$.} \endproclaim

\noindent{\it Proof.} Let $\bar M(E_\varrho) :=
M(E_\varrho)\cup\{\infty\}$ denote the one point compactification
of $M(E_\varrho)$. Then the sequence $\{F_n\}$ viewed as
probabilities on $\bar M(E_\varrho)$ is relatively compact. Let
$\{F_{n_k}\}$ be a subsequence of $\{F_n\}$ which converges to
some probability measure $F$ on $\bar M(E_\varrho)$. By
(\ref{2.7}) and bounded convergence, we have
 \beqnn
L(f) = \int_{\bar M(E_\varrho)}\e^{-\nu(f)}F(\d\nu), \quad f\in
D(E_\varrho)^+,
 \eeqnn
where the integrand is defined as zero at $\infty$ by continuity.
Since $L(f)\to 1$ as $f\to 0$, we have $F(M(E_\varrho)) =1$ and
hence the first assertion follows. The second assertion is
immediate.\qed

For any $\alpha\in [-\infty,\infty)$, a family of $\sigma$-finite
measures $(K_t)_{t>\alpha}$ on $M(E)$ is called an {\it entrance
law} (at $\alpha$) for the semigroup $(Q_t)_{r\le t}$ if $K_r
Q_{r,t} =K_t$ for all $t>r>\alpha$. It is called a {\it
probability entrance law} if each $K_t$ is a probability measure.
An entrance law $(K_t)_{t>\alpha}$ is said to be {\it closable} if
there is a $\sigma$-finite measure $K_\alpha$ on $M(E)$ such that
$K_t = K_\alpha Q_{\alpha,t}$ for all $t>\alpha$. In this case,
$(K_t)_{t\ge \alpha}$ is called a {\it closed entrance law} for
$(Q_{r,t})_{r\le t}$. An entrance law $(K_t)_{t>\alpha}$ is said
to be {\it minimal} if every entrance law dominated by
$(K_t)_{t>\alpha}$ is proportional to it. Those definitions are
applicable to general transition semigroups under obvious
modifications.

\medskip
\noindent{\bf Example 2.1.} Let $T_1\subset \IR$ be a countable
set and $\{(K_{s,t})_{t>s}: s\in T_1\}$ be a family of probability
entrance laws for $(Q_{r,t})_{r\le t}$. Suppose that
 \beqnn
-\sum_{s\in[r,t)\cap T_1}\log L_{K_{s,t}}(1) <\infty
 \eeqnn
for all $r\le t\in \IR$. By Lemma\n{2.2} we may see that
 \beqlb
\log L_{N_{r,t}}(f)
=
\sum_{s\in[r,t)\cap T_1}\log L_{K_{s,t}}(f),
\quad r\le t,\ f\in B(E)^+,
\label{2.8}
 \eeqlb
defines a probability measure $N_{r,t}$ on $M(E)$. It is simple to
check that $(N_{r,t}) _{r\le t}$ form an SC-semigroup. Suppose
that for each $s\in T_1$ we have a Markov process $(X_{s,t})_{t>
s}$ with transition semigroup $(Q_{r,t})_{r\le t}$ and
one-dimensional distributions $(K_{s,t})_{t>s}$ and that the
family $\{(X_{s,t})_{t> s}: s\in T_1\}$ are independent. Then for
any $t\ge \alpha$ the random measure
 \beqnn
Y_t := \sum_{s\in[a,t)\cap T_1} X_{s,t}
 \eeqnn
is a.s. well-defined and $\{Y_t: t\ge \alpha\}$ is an immigration
process corresponding to the SC-semigroup given by (\ref{2.8}).

\medskip
\noindent{\bf Example 2.2.} Let $T_2\subset \IR$ be a countable
set and $\{(K_{s,t})_{t\ge s}: s\in T_2\}$ be a family of closed
probability entrance laws for $(Q_{r,t}) _{r\le t}$ such that
 \beqnn
-\sum_{s\in(r,t]\cap T_2}\log L_{K_{s,t}}(1) <\infty
 \eeqnn
for all $r\le t \in \IR$. Then we may define an SC-semigroup
$(N_{r,t})_{r\le t}$ by
 \beqlb
\log L_{N_{r,t}}(f)
=
\sum_{s\in(r,t]\cap T_2}\log L_{K_{s,t}}(f),
\quad r\le t,\ f\in B(E)^+.
\label{2.9}
 \eeqlb
The corresponding immigration process can be constructed similarly
as in the last example.

\medskip
\noindent{\bf Example 2.3.} Suppose that $(N_{r,t})_{r\le t}$ is a
family of probability measures on $M(E)$ given by
 \beqlb
\log L_{N_{r,t}}(f) = \int_r^t\log L_{K_{s,t}}(f)\zeta(\d s),
\quad r\le t,\ f\in B(E)^+, \label{2.10}
 \eeqlb
where $\zeta(\d s)$ is a Radon measure on $\IR$ and
$\{(K_{s,t})_{t>s}: s\in\IR\}$ is a family of infinitely divisible
probability entrance laws for $(Q_{r,t})_{r\le t}$. Then
$(N_{r,t})_{r\le t}$ is an SC-semigroup.

\medskip
\noindent{\bf Example 2.4.} Suppose we have a probability space on
which the two processes $\{X^1_t: t\ge 0\}$ and $\{X^2_t: t\ge
\alpha\}$ are defined, where $\{X^1_t: t\ge 0\}$ is a superprocess
with parameters $(\xi_1,\phi_1)$, and $\{X^2_t: t\ge 0\}$
conditioned on $\{X^1_t: t\ge 0\}$ is an immigration process
corresponding to the SC-semigroup determined by
 \beqnn
\log L_{N_{r,t}}(f)
=
- \int_r^tX_s(V^2_{t-s}f)\d s,
\quad r\le t,\ f\in B(E)^+,
 \eeqnn
where $(V^2_t)_{t\ge 0}$ is defined by (\ref{1.2}) with parameters
$(\xi_2,\phi_2)$. Then $\{(X^1_t, X^2_t): t\ge 0\}$ is also a
Markov process. Intuitively, it describes the evolution of a
population with two types of ``particles'' on $E$, where the
second type can be produced by both of them; see Hong and Li
\cite{HL99}. More general forms of multi-type superprocesses were
studied in Gorostiza and Lopez-Mimbela \cite{GL90}, Gorostiza and
Roelly \cite{GR90}, Li \cite{L92a}, etc.

\section{Decomposition of skew convolution semigroups}

\setcounter{equation}{0}

In this section, we prove a decomposition theorem for the
inhomogeneous SC-semigroup, which appears even in the construction
of homogeneous immigration processes. Let us consider the
transition semigroup $(Q_{r,t})_{r\le t}$ of an MB-process defined
in the last section.

\proclaim{\bf Theorem\n{3.1}.} {\sl A family of probability
measures $(N_{r,t})_{r\le t}$ on $M(E)$ form an SC-semigroup
associated with $(Q_{r,t})_{r\le t}$ if and only if
 \beqlb
\log L_{N_{r,t}}(f)
&=&
\sum_{s\in[r,t)\cap T_1}\log L_{K^1_{s,t}}(f)
+ \sum_{s\in(r,t]\cap T_2}\log L_{K^2_{s,t}}(f) \nonumber \\
& &\qquad + \int_r^t\log L_{K^3_{s,t}}(f)\zeta(\d s), \qquad r\le
t,\ f\in B(E)^+, \label{3.1}
 \eeqlb
where $T_1, T_2\subset \IR$ are countable sets, $\zeta(\d s)$ is a
diffuse Radon measure on $\IR$, $\{(K^1_{s,t})_{t>s}: s\in T_1\}$
is a family of probability entrance laws, $\{(K^2_{s,t})_{t\ge s}:
s\in T_2\}$ is a family of closed probability entrance laws, and
$\{(K^3_{s,t})_{t>s}: s\in\IR\}$ is a family of infinitely
divisible probability entrance laws for $(Q_{r,t})_{r\le t}$.}
\endproclaim

In principle, the entrance laws are obtained via applications of
Lemma\n{2.2}. The proof is a little tedious because we are not
assuming the Feller property and the class $D(E_\varrho)^+$ is not
preserved by the cumulant semigroup $(V_{r,t})_{r\le t}$. We shall
break the proof into several lemmas. Suppose that $(N_{r,t})_{r\le
t}$ is an SC-semigroup associated with $(Q_{r,t})_{r\le t}$ and
let
 \beqlb
J_{r,t}(f) = -\log \int_{M(E)} \e^{-\nu(f)} N_{r,t}(\d\nu),
\quad r\le t,\ f\in B(E)^+.
\label{3.2}
 \eeqlb
Then the relation (\ref{2.5}) is equivalent to
 \beqlb
J_{r,t}(f) = J_{r,s}(V_{s,t}f) + J_{s,t}(f),
\quad r\le s\le t,\ f\in B(E)^+.
\label{3.3}
 \eeqlb
By (\ref{3.3}) one sees that $J_{r,t}(f)$ is a non-increasing
function of $r\le t$ and $J_{t,t}(f)$ =0. Then there is a Radon
measure $G_t(f,\cdot)$ on $(-\infty,t]$ such that $G_t(f,\{t\}) =
\lim_{r\uparrow t}J_{r,t}(f)$ and
 \beqlb
G_t(f,(r,s]) = \lim_{v\downarrow s}\lim_{u\downarrow r}[J_{u,t}(f)
- J_{v,t}(f)] = \lim_{v\downarrow s}\lim_{u\downarrow
r}J_{u,v}(V_{v,t}f), \qquad r\le s< t. \label{3.4}
 \eeqlb

\proclaim{\bf Lemma\n{3.2}.} {\sl (i)~If $f,g\in B(E)^+$ and $0\le
g\le cf$ for a constant $c\ge 1$, then $G_t(g,\cdot) \le
cG_t(f,\cdot)$. (ii)~For $t\le u$ and $f\in B(E)^+$, we have
$G_u(f,\d s) = G_t(V_{t,u}f,\d s)$ on $(-\infty,r)$. (iii)~For
$t\le u$ and $f\in B(E)^+$, we have $G_u(f,\d s) \ll G_t(1,\d s)$
on $(-\infty,t)$.} \endproclaim

\noindent{\it Proof.} For $r\le s< t\le u$, we use (\ref{2.4}),
(\ref{3.2}) and Jensen's inequality to see that
 \beqnn
G_t(g,(r,s]) &=& \lim_{v\downarrow s}\lim_{w\downarrow
r}J_{w,v}(V_{v,t}g)
\le \lim_{v\downarrow s}\lim_{w\downarrow r}J_{w,v}(V_{v,t}(cf)) \\
&\le& c\lim_{v\downarrow s}\lim_{w\downarrow r}J_{w,v}(V_{v,t}f) =
cG_t(f,(r,s]).
 \eeqnn
Similarly, we have $G_t(g,\{t\}) \le cG_t(f,\{t\})$. Then (i)
follows. Let $r\le s< t\le u$. By (\ref{3.4}) and the semigroup
property of $(V_{r,t})_{r\le t}$,
 \beqnn
G_u(f,(r,s]) = \lim_{v\downarrow s}\lim_{w\downarrow
r}J_{w,v}(V_{v,t}V_{t,u}f) = G_t(V_{t,u}f,(r,s]),
 \eeqnn
yielding (ii). Combining (i) and (ii) we have
 \beqnn
G_u(f,(r,s])
= G_t(V_{t,u}f,(r,s])
\le (\|V_{t,u}f\|+1)G_t(1,(r,s]),
 \eeqnn
from which (iii) follows. \qed

Clearly, we have the unique decomposition:
 \beqlb
J_{r,t}(f) = J_{r,t}^1(f) + J_{r,t}^2(f) + J_{r,t}^3(f),
\quad s\le t,\ f\in B(E)^+,
\label{3.5}
 \eeqlb
where $J_{r,t}^1(f)$ is a left continuous non-increasing step
function, $J^2_{r,t} (f)$ is a right continuous non-increasing
step function and $J_{r,t}^3(f)$ is a continuous non-increasing
function of $r\le t$, and $J_{t,t}^i(f)=0$ for $i=1,2,3$. By the
uniqueness, $(J_{r,t}^i)_{r\le t}$ also satisfies equation
(\ref{3.3}). Applying Lemma\n{2.2} we may get
 \beqlb
N_{r,t}
=
N^1_{r,t} *N^2_{r,t} *N^3_{r,t},
\quad r\le t,
\label{3.6}
 \eeqlb
where $(N^i_{r,t}) _{r\le t}$ is the SC-semigroup corresponding to
the functional $(J_{r,t}^i)_{r\le t}$. Observe also that
$J_{s,t}^1(f)$ and $J_{s,t}^2(f)$ yield atomic measures
$G_t^1(f,\cdot)$ and $G_t^2(f,\cdot)$ on $(-\infty,t]$,
respectively, and $J_{s,t}^3(f)$ yields a diffuse measure
$G_t^3(f,\cdot)$ on $(-\infty,t]$.

\proclaim{\bf Lemma\n{3.3}.} {\sl There are countable sets $T_1,
T_2\subset\IR$, probability entrance laws $\{(K^1_{r,t}) _{t>s}:
r\in T_1\}$ and closed probability entrance laws $\{(K^2_{r,t})
_{t\ge r}: r\in T_2\}$ such that
 \beqlb
J^1_{r,t}(f)
=
-\sum_{s\in[r,t)\cap T_1}\log L_{K^1_{s,t}}(f),
\quad r\le t,\ f\in B(E)^+,
\label{3.7}
 \eeqlb
and
 \beqlb
J^2_{r,t}(f)
=
-\sum_{s\in(r,t]\cap T_2}\log L_{K^2_{s,t}}(f),
\quad r\le t,\ f\in B(E)^+.
\label{3.8}
 \eeqlb
}\endproclaim

\noindent{\it Proof.} Since the arguments are similar, we only
give the proof of (\ref{3.8}). Recall that $J^2_{s,t}(f)$ is right
continuous in $s\le t$. Consequently,
 \beqlb
G_t^2(f,\{s\}) = \downarrow\lim_{r\uparrow s}[J^2_{r,t}(f) -
J^2_{s,t}(f)] = \downarrow\lim_{r\uparrow s}J^2_{r,s}(V_{s,t}f),
\quad f\in B(E)^+. \label{3.9}
 \eeqlb
That is,
 \beqnn
\exp\{-G_t^2(f,\{s\})\} = \uparrow\lim_{r\uparrow
s}\int_{M(E)}\e^{-\nu(f)}N^2_{r,s}Q_{s,t}(\d \nu), \quad f\in
B(E)^+.
 \eeqnn
By Lemma\n{2.2} we see that
 \beqlb
\exp\{-G_t^2(f,\{s\})\} = \int_{M(E)}\e^{-\nu(f)}K^2_{s,t}(\d
\nu), \quad f\in D(E_\varrho)^+, \label{3.10}
 \eeqlb
for a probability measure $K^2_{s,t}$ on $M(E)$. Let $Q= \{u_i:
i=1,2,\cdots\}$ be a countable dense subset of $\IR$ and let $T_2$
be the collection of atoms of the measures $\{G^2_{u_i} (1,\cdot):
i=1,2,\cdots\}$. By (\ref{3.10}) and Lemma\n{3.2},
 \beqnn
\log L_{N^2_{r,t}}(f)
=
- G_t^2(f,(r,t])
=
\log L_{K^2_{t,t}}(f)
+ \sum_{s\in(r,t)\cap T_2}\log L_{K^2_{s,t}}(f),
\quad r\le t,
 \eeqnn
first for $f\in D(E_\varrho)^+$ and then for all $f\in B(E)^+$. In
particular, (\ref{3.10}) also holds for $f\in B(E)^+$. But,
(\ref{3.9}) implies that $G_t^2(f,\{s\}) = G_s^2 (V_{s,t}f,
\{s\})$. Then we must have $K_{s,s}Q_{s,t} = K_{s,t}$ for $s\le
t$, that is $(K_{s,t})_{t\ge s}$ form a closed entrance law. If
$-\log L_{K^2_{t,t}}(f) >0$, then $-\log L_{K^2_{t,t}}(1) >0$. By
the continuity of $V_{t,u}1$ in $u\ge t$,
 \beqnn
G_u^2(1,\{t\})
=
- \log L_{K^2_{t,u}}(1)
=
- \log L_{K^2_{t,t}}(V_{t,u}1)
>0
 \eeqnn
for some $u\in Q$. Then we have $t\in Q$ and (\ref{3.8}) follows.
\qed

\proclaim{\bf Lemma\n{3.4}.} {\sl There is a diffuse Radon measure
$\zeta(\d s)$ on $\IR$ such that $G_t^3(f,\d s)\ll \zeta(\d s)$
for all $t\in\IR$ and $f\in B(E)^+$.}
\endproclaim

\noindent{\it Proof.} Let $Q= \{u_i: i=1,2,\cdots\}$ be a
countable dense subset of $\IR$ and choose $a_i>0$ such that
 \beqnn
\sum_{i=1}^\infty a_iG^3_{u_i}(1,(u_i-1,u_i]) <\infty.
 \eeqnn
We may define a diffuse Radon (in fact finite) measure $\zeta$ on
$\IR$ by
 \beqlb
\zeta(\d s) = \sum_{i=1}^\infty 1_{(u_i-1,u_i]}(s)G^3_{u_i}(1,\d
s), \quad s\in \IR. \label{3.11}
 \eeqlb
If $G^3_t(f,B) >0$ for a Borel set $B\subset (-\infty,t]$, then
$G^3_t(f, (u_i-1,u_i)\cap B) >0$ for some $u_i\in (-\infty,t)\cap
Q$. Therefore, Lemma\n{3.2} implies that $G^3_{u_i} (1,
(u_i-1,u_i)\cap B) >0$, and hence $\zeta(B)>0$. \qed

\proclaim{\bf Lemma\n{3.5}.} {\sl The Radon-Nikodym derivative
$G^3_t(f,\d s)/\zeta(\d s)$ has a version $D_{s,t}(f)$ with the
representation
 \beqlb
D_{s,t}(f)
=
-\log L_{F_{s,t}}(f),
\quad s<t,\ f\in B(E)^+,
\label{3.12}
 \eeqlb
where $F_{s,t}$ is an infinitely divisible probability measure on
$M(E)$.}
\endproclaim

\noindent{\it Proof.} Let $r<t$ and assume $\zeta(r,t] >0$ to
avoid triviality. We take an increasing sequence of ordered sets
$\pi_n = \{s^n_0, s^n_1, \cdots, s^n_{m(n)}\}$ with $r= s^n_0<
s^n_0<\cdots<s^n_{m(n)} =t$ and $\lim_{n\to\infty}\delta_n =0$,
where $\delta_n = \max \{s^n_i -s^n_{i-1}: 1\le i\le m(n)\}$. Let
${\cal F}_n$ be the $\sigma$-algebra on $(r,t]$ generated by
$(s^n_0, s^n_1], \cdots, (s^n_{m(n)-1},s^n_{m(n)}]$ and let
 \beqnn
M_n(f)(s) =\left\{\begin{array}{ll}
{G^3_t(f,(s^n_{i-1},s^n_i])}/{\zeta(s^n_{i-1},s^n_i]}
&\mbox{ if }s^n_{i-1}<s\le s^n_i \mbox{ and }\zeta(s^n_{i-1},s^n_i] >0,  \\
0 &\mbox{ if }s^n_{i-1}<s\le s^n_i \mbox{ and
}\zeta(s^n_{i-1},s^n_i] =0.
\end{array}\right.
 \eeqnn
Then $\{M_n(f),{\cal F}_n: n\ge 1\}$ under the probability measure
$\zeta(r,t] ^{-1}\zeta$ is a martingale which is closed on the
right by the Radon-Nikodym derivative $G^3_t(f,\d s) /\zeta(\d
s)$. But, since $\{{\cal F}_n: n\ge 1\}$ generates the Borel
$\sigma$-algebra on $(r,t]$, $M_n(f)(s)$ converges as $n\to\infty$
to $G^3_t(f,\d s)/\zeta(\d s)$ for $\zeta$-a.e. $s\in (r,t]$ by
the martingale convergence theorem. Then we may find a set
$A_t\subset (-\infty,t]$ with full $\zeta$-measure such that for
any $s\in A_t$ there are sequences $r_k\uparrow s$ and
$s_k\downarrow s$ satisfying
 \beqnn
G^3_t(f,\d s)/\zeta(\d s) =
\lim_{k\to\infty}G^3_t(f,(r_k,s_k])/\zeta(r_k,s_k], \quad f\in
D(E_\varrho)^+.
 \eeqnn
Clearly, $G^3_t(f,(r_k,s_k]) = G_{s_k}(V_{s_k,t}f,(r_k,s_k])$ goes to zero as
$k\to\infty$. It follows that
 \beqlb
\frac{G^3_t(f,\d s)}{\zeta(\d s)} &=& \lim_{k\to\infty}
\frac{1-\exp\left\{-G_{s_k}(V_{s_k,t}f,(r_k,s_k])\right\}}{\zeta(r_k,s_k]}
\nonumber \\
&=& \lim_{k\to\infty}\frac{1}{\zeta(r_k,s_k]}\int_{M(E)^\circ}
\left(1-\e^{-\nu(f)}\right) N_{r_k,s_k}Q_{s_k,t}(\d\nu)
\label{3.13}
 \eeqlb
for $f\in D(E_\varrho)^+$. On the other hand, since
 \beqnn
\int_r^t[G^3_t(f,\d s)/\zeta(\d s)]\zeta(\d s) = - \log
\int_{M(E)}\e^{-\nu(f)}N^3_{r,t}(\d\nu),
 \eeqnn
choosing a smaller full $\zeta$-measure set $A_t\subset
(-\infty,t]$ we may get $G^3_t(f,\d s)/\zeta(\d s)\downarrow0$ as
$f\downarrow0$ for all $s\in A_t$. By (\ref{3.13}) and
Lemma\n{2.2} there are infinitely divisible probability measures
$\{F_{s,t}: s\in A_t\}$ such that
 \beqnn
G^3_t(f,\d s)/\zeta(\d s) = -\log L_{F_{s,t}}(f), \quad s\in A_t,\
f\in D(E_\varrho)^+.
 \eeqnn
Setting $F_{s,t} = \delta_0$ for $s\in (-\infty,t]\setminus A_t$
we have
 \beqnn
G^3_t(f,(r,t]) = - \log \int_{M(E)}\e^{-\nu(f)}N^3_{r,t}(\d\nu) =
-\int_r^t\log L_{F_{s,t}}(f)\zeta(\d s),
 \eeqnn
first for $f\in D(E_\varrho)^+$ and then for all $f\in B(E)^+$.
\qed

\proclaim{\bf Lemma\n{3.6}.} {\sl There are infinitely divisible
probability entrance laws $\{(K_{s,t})_{t>s}: s\in\IR\}$ for the
semigroup $(Q_{r,t})_{r\le t}$ such that
 \beqlb
J_{r,t}^3(f) = -\int_r^t \log L_{K_{s,t}}(f) \zeta(\d s), \quad
t\ge r,\ f\in B(E)^+. \label{3.14}
 \eeqlb
} \endproclaim

\noindent{\it Proof.} By Lemma\n{3.2} we have $G^3_t(f,\d s) =
G^3_r (V_{r,t}f,\d s)$ for $s< r\le t$. It follows that
 \beqnn
D_{s,t}(f) = D_{s,r}(V_{r,t}f), \quad \zeta\mbox{-a.e.
}s\in(-\infty,r].
 \eeqnn
Let $\lambda$ denote the Lebesgue measure on $\IR$. By Fubini's
lemma, there is a set $B(f)\subseteq \IR$ with full
$\zeta$-measure and sets $C_s(f)\subseteq [s,\infty)$ and
$C_{s,r}(f)\subseteq [r,\infty)$ with full $\lambda$-measure such
that
 \beqnn
D_{s,t}(f) = D_{s,r}(V_{r,t}f),
\quad s\in B(f),\ r\in C_s(f),\ t\in C_{s,r}(f).
 \eeqnn
By Lemma\n{3.5}, $D_{s,t}$ and $D_{s,r}\circ V_{r,t}$ are
determined by their restrictions to the countable class
$D(E_\varrho)^+$. Then for a set $B\subseteq \IR$ with full
$\zeta$-measure and sets $C_s\subseteq (s,\infty)$ and
$C_{s,r}\subseteq (r,\infty)$ with full $\lambda$-measures we have
 \beqlb
D_{s,t} = D_{s,r}\circ V_{r,t},
\quad s\in B,\ r\in C_s,\ t\in C_{s,r},
\label{3.15}
 \eeqlb
as operators on $B(E)^+$. For any $s\in B$, choose a sequence
$\{s_k\} \subset C_s$ with $s_k\downarrow s$. By (\ref{3.15}) we
get
 \beqnn
D_{s,s_k} \circ V_{s_k,t}
= D_{s,s_{k+1}} \circ V_{s_{k+1},t}
= D_{s,t},
\quad t\in C_{s,s_k}\cap C_{s,s_{k+1}}.
 \eeqnn
By Lemma\n{3.5} and our assumptions, $D_{s,s_k} (V_{s_k,t}f)$ and
$D_{s,s_{k+1}} (V_{s_{k+1},t}f)$ are right continuous in $t\ge
s_{k+1}$ for $f\in C(E)^+$. Then we have
 \beqlb
D_{s,s_k} \circ V_{s_k,t} = D_{s,s_{k+1}} \circ V_{s_{k+1},t},
\quad t\ge s_k,
\label{3.16}
 \eeqlb
as operators on $B(E)^+$. For $t>s$ we take some $s_k\in (s,t)$
and let $I_{s,t} = D_{s,s_k} \circ V_{s_k,t}$, which is
independent of the choice of $s_k$ by (\ref{3.16}).
Correspondingly, let $F_{s,s_k}$ be the infinitely divisible
probability measure on $M(E)$ given by Lemma\n{3.5} and let
$K_{s,t} = F_{s,s_k}Q_{s_k,t}$, which is independent of $s_k$
ether. Clearly, $(K_{s,t})_{t>s}$ form an entrance law for
$(Q_{r,t})_{r\le t}$ and
 \beqlb
I_{s,t} = -\log L_{F_{s,s_k}}\circ V_{s_k,t}
=
-\log L_{K_{s,t}}(f),
\quad t>s.
\label{3.17}
 \eeqlb
From (\ref{3.15}) we have $I_{s,t} = D_{s,t}$ for $\lambda$-a.e.
$t\in (s,\infty)$. If $s\in \IR\setminus B$, we define $I_{s,t}=0$
for all $t>s$. Then Fubini's lemma yields the existence of a set
$U\subseteq \IR$ with full $\lambda$-measure such that for each
$t\in U$ we have $I_{s,t} = D_{s,t}$ for $\zeta$-a.e.
$s\in(-\infty,t)$. Consequently, if $t\in U$,
 \beqlb
J_{r,t}^3(f) = G^3_t(f,(r,t]) = \int_r^t I_{s,t}(f) \zeta(\d s),
\quad r\le t,\ f\in B(E)^+. \label{3.18}
 \eeqlb
Recall that both $\zeta$ and $G^3_t(f,\cdot)$ are diffuse
measures. For $t\in\IR\setminus U$ we choose a sequence
$\{t_k\}\subset U$ with $t_k\uparrow t$. By (\ref{3.18}) and
Lemma\n{3.2},
 \beqnn
J_{r,t}^3(f)
&=& \lim_{k\to\infty}G^3_t(f,(r,t_k])
= \lim_{k\to\infty}G^3_{t_k}(V_{t_k,t}f,(r,t_k]) \\
&=& \lim_{k\to\infty}\int_r^{t_k} I_{s,t_k}(V_{t_k,t}f) \zeta(\d
s)
= \lim_{k\to\infty}\int_r^{t_k} I_{s,t}(f) \zeta(\d s) \\
&=& \int_r^t I_{s,t}(f) \zeta(\d s),
 \eeqnn
yielding the desired result.
\qed

\noindent{\it Proof of Theorem\n{3.1}.} It is easy to see that, if
$(N_{r,t})_{r\le t}$ is given by (\ref{3.1}), it satisfies
(\ref{2.5}). Combining (\ref{3.5}) and Lemmas\n{3.3} and\n{3.6} we
see that any SC-semigroup has the decomposition (\ref{3.1}). \qed

Let $(Q_t)_{t\ge 0}$ be the transition semigroup of a homogeneous
MB-process. The following special case of Theorem\n{3.1} was
proved in Li \cite{L96a}.

\proclaim{\bf Theorem\n{3.7} {\rm (Li \cite{L96a})}.} {\sl A
family of probability measures $(N_t)_{t\ge0}$ on $M(E)$ is an
SC-semigroup associated with $(Q_t)_{t\ge0}$ if and only if there
is an infinitely divisible probability entrance law $(K_t)_{t>0}$
for $(Q_t)_{t\ge0}$ such that
 \beqlb
\log L_{N_t}(f)
=
\int_0^t\big[\log L_{K_s}(f)\big]\d s,
\quad t\ge 0,\ f\in B(E)^+.
\label{3.19}
 \eeqlb
} \endproclaim
\qed

\section{Construction of immigration processes}
\setcounter{equation}{0} To make best use of the existing
literature, we restrict in this and the subsequent sections to the
transition semigroup $(Q_t)_{t\ge0}$ of a homogeneous Borel right
MB-process $X$. Let $(Q_t^\circ)_{t\ge0}$ denote the restriction
of $(Q_t)_{t\ge0}$ to $M(E)^\circ$.

We first review some facts in potential theory; see e.g.
Dellacherie et al \cite{DMM92} and Getoor \cite{G90}. A family of
$\sigma$-finite measures $(J_t)_{t\in\IR}$ is called an {\it
entrance rule} for $(Q_t^\circ)_{t\ge0}$ if $J_sQ_{t-s}^\circ\le
J_t$ for $t>s\in\IR$ and $J_s Q_{t-s}^\circ\uparrow J_t$ as
$s\uparrow t$. Note that an entrance law $(H_t)_{t>r}$ at $r\in
[-\infty,\infty)$ can be extended to an entrance rule by setting
$H_t=0$ for $t\le r$. Let $W(M(E))$ denote the space of paths
$\{w_t: t\in\IR\}$ that are $M(E)^\circ$-valued and right
continuous on an open interval $(\alpha(w),\beta(w))$ and take the
value of the null measure elsewhere. The path $[0]$ constantly
equal to $0$ corresponds to $(\alpha,\beta)$ being empty. Set
$\alpha([0]) = +\infty$ and $\beta([0]) = -\infty$. Let $({\cal
G}^\circ, {\cal G}_t^\circ)_{t\in\IR}$ be the natural
$\sigma$-algebras on $W(M(E))$ generated by the coordinate
process. The shift operators $\{\sigma_t: t\in\IR\}$ on $W(M(E))$
are defined by $\sigma_tw_s=w_{t+s}$. To any entrance rule
$(J_t)_{t\in\IR}$ there corresponds a unique $\sigma$-finite
measure $\bQ^J$ on $(W(M(E)),{\cal G}^\circ)$ under which the
coordinate process $\{w_t: t\in\IR\}$ is a Markov process with
one-dimensional distributions $(J_t)_{t\in\IR}$ and semigroup
$(Q_t^\circ)_{t\ge0}$. That is,
 \beqlb
&&\bQ^J\{\alpha<t_1,w_{t_1}\in\d\nu_1, w_{t_2}\in\d\nu_2, \cdots,
w_{t_n}\in\d\nu_n, t_n<\beta\} \nonumber \\
&&\qquad = J_{t_1}(\d\nu_1)Q_{t_2-t_1}^\circ(\nu_1,\d\nu_2)\cdots
Q_{t_n-t_{n-1}}^\circ(\nu_{n-1},\d\nu_n) \label{4.1}
 \eeqlb
for all $t_1 < \cdots < t_n \in \IR$ and $\nu_1, \cdots, \nu_n\in
M(E)^\circ$. The existence of this measure was proved by Kuznetsov
\cite{K74}; see also Getoor and Glover \cite{GG87}. The system
$(W(M(E)), {\cal G}^\circ, {\cal G}^\circ_t$, $w_t,\bQ^J)$ is now
commonly called the {\it Kuznetsov process} determined by
$(J_t)_{t\in\IR}$, and $\bQ^J$ is called the {\it Kuznetsov
measure}. We have the representation
 \beqlb
J_t = H_{-\infty,t}+\int_{\IR} H_{s,t}\rho(\d s),
\quad t\in\IR,
\label{4.2}
 \eeqlb
where $\rho(\d s)$ is a Radon measure on $\IR$ and
$(H_{s,t})_{t\in\IR}$ is an entrance law at $s\in
[-\infty,\infty)$. This representation yields
 \beqlb
\bQ^J(\d w) = \bQ_{-\infty}(\d w) + \int_{\IR}\bQ_s(\d w)\rho(\d
s), \label{4.3}
 \eeqlb
where $\bQ_s(\d w)$ is the Kuznetsov measure determined by
$(H_{s,t})_{t\in\IR}$; see \cite{GG87}. If $(J_t)_{t\in\IR}$ is an
entrance law at $r\in [-\infty,\infty)$, then $\bQ^J$ is supported
by $W_r(M(E))$, the subset of $W(M(E))$ comprising paths $\{w_t:
t\in\IR\}$ such that $\alpha(w)=r$. In particular, if $F$ is an
excessive measure for $(Q^\circ_t)_{t\ge0}$ and $J_t \equiv F$,
then $\bQ^J$ is stationary, that is, $\bQ^J\circ\sigma_t^{-1} =
\bQ^J$ for all $t\in\IR$.

Now suppose that $(J_t)_{t\in\IR}$ is an entrance rule for
$(Q_t^\circ)_{t\ge0}$ with the representation (\ref{4.2}) and
$N^J(\d w)$ is a Poisson random measure on $W(M(E))$ with
intensity $\bQ^J(\d w)$. It is easy to see that
 \beqlb
Y_t^J := \int_{W(M(E))} w_t\, N^J(\d w)
\label{4.4}
 \eeqlb
is a.s. well-defined for each $t\in\IR$.

\proclaim{\bf Lemma\n{4.1}.} {\sl In the situation described
above, $\{Y^J_t: t\in\IR\}$ is an immigration process
corresponding to the (inhomogeneous) SC-semigroup $(N_{r,t})_{r\le
t}$ defined by
 \beqlb
\log L_{N_{r,t}}(f) = - \int_{[r,t)}\int_{M(E)^\circ}
\left(1-\e^{-\nu(f)}\right) H_{s,t}(\d\nu)\rho(\d s), \quad r\le
t, f\in B(E)^+. \label{4.5}
 \eeqlb
} \endproclaim

\noindent{\it Proof.} By (\ref{4.3}), for any bounded Borel
function $F$ on $M(E)$ with $F(0)=0$, we have
 \beqlb
\bQ^J\big\{F(w_t); r\le \alpha <t\big\} =
\int_{[r,t)}H_{s,t}(F)\,\rho(\d s). \label{4.6}
 \eeqlb
Then the results follow from (\ref{4.6}) and the Markov property
of $\bQ^J$. \qed

Let $(N_t)_{t\ge0}$ an SC-semigroup associated with $(Q_t)_{t\ge
0}$ which is given by (\ref{3.19}). Suppose that $K_t=
I(\eta_t,H_t)$ and $N_t= I(\gamma_t,G_t)$ for $t>0$ and $t\ge0$,
respectively.

\proclaim{\bf Lemma\n{4.2}.} {\sl Let $G_t=0$ for $t<0$. Then
$(G_t)_{t\in\IR}$ is an entrance rule for $(Q_t^\circ)_{t\ge0}$.}
\endproclaim

\noindent{\it Proof.} Recall that $Q_t(\mu,\cdot)$ is an
infinitely divisible probability measure on $M(E)$ for all $t\ge
0$ and $\mu\in M(E)$. Suppose $Q_t(\delta_x,\cdot) =
I(\lambda_t(x,\cdot), L_t(x,\cdot))$. By (\ref{1.4}),
 \beqnn
G_t = G_{t-r} + G_rQ_{t-r}^\circ + \int_E\gamma_r(\d
x)L_{t-r}(x,\cdot), \quad t>r>0,
 \eeqnn
and hence $G_rQ_{t-r}\le G_t$. From (\ref{3.19}) we have
 \beqnn
G_tQ_{t-r}^\circ
= \int_0^tH_sQ_{t-r}^\circ\d s
= G_rQ_{t-r}^\circ + \int_r^tH_sQ_{t-r}^\circ\d s,
 \eeqnn
so $G_rQ_{t-r}^\circ\uparrow G_t$ as $r\uparrow t$. Therefore,
$(G_t)_{t\in\IR}$ is an entrance rule. \qed

Now we give the construction of the immigration process
corresponding to $(N_t)_{t\ge0}$. The next theorem shows that,
except the deterministic part $\{\gamma_t: t\ge 0\}$, both the
entering times and the evolutions of the immigrants are decided by
a Poisson random measure based on the Kuznetsov measure $\bQ^G$.

\proclaim{\bf Theorem\n{4.3}.} {\sl Let $Y^G_t$ be defined by
(\ref{4.4}) with $J=G$ and let $Y_t = \gamma_t + Y^G_t$. Then
$\{Y_t: t\ge0\}$ is an immigration process with one-dimensional
distributions $(N_t)_{t\ge0}$ and transition semigroup
$(Q^N_t)_{t\ge0}$.} \endproclaim

\noindent{\it Proof.} Suppose $(G_t)_{t\in\IR}$ is represented by
(\ref{4.2}) with $G_{s,t}$ in place of $J_{s,t}$. Then, for $t\ge
r\ge 0$,
 \beqlb
\int_{[r,t)}G_{s,t}\,\rho(\d s)
&=&
\int_{[0,t)}G_{s,t}\,\rho(\d s)
- \int_{[0,r)}G_{r,s} Q_{t-r}^\circ\rho(\d s) \nonumber \\
&=&
G_t - G_rQ_{t-r}^\circ \nonumber \\
&=&
\int_0^tH_s \d s - \int_0^rH_sQ_{t-r}^\circ\d s.
\label{4.7}
 \eeqlb
The relation $K_{s+t}=K_sQ_t$ yields
 \beqlb
\eta_{s+t} = \int_E\eta_s(\d x)\lambda_t(x,\cdot), \quad H_{s+t} =
\int_E \eta_s(\d x)L_t(x,\cdot) + H_sQ_t^\circ. \label{4.8}
 \eeqlb
From the second equation in (\ref{4.8}) we have
 \beqnn
\int_0^r H_{s+t-r} \d s - \int_0^rH_sQ_{t-r}^\circ\d s =
\int_0^r\d s\int_E \eta_s(\d x) L_{t-r}(x,\cdot).
 \eeqnn
Substituting this into (\ref{4.7}) gives
 \beqlb
\int_{[r,t)}G_t^s\,\rho(\d s) &=& \int_0^{t-r}H_s \d s +
\int_0^r\d s\int_E \eta_s(\d x)
L_{t-r}(x,\cdot) \nonumber \\
&=& G_{t-r} + \int_E \gamma_r(\d x) L_{t-r}(x,\cdot). \label{4.9}
 \eeqlb
Since $\{\gamma_t: t\ge0\}$ is deterministic, it is simple to
check that $\{Y_t: t\ge0\}$ is a Markov process with
one-dimensional distributions $(N_t)_{t\ge0}$. By Lemma\n{4.1} we
have
 \beqlb
&&\bE [\exp\{- Y_t(f)\}\big| Y_s: 0\le s\le r] \nonumber \\
&&\qquad
= \exp\bigg\{-Y^G_r(V_{t-r}f)-\gamma_t(f) \nonumber \\
&&\qquad\qquad\qquad - \int_{[r,t)}\rho(\d
s)\int_{M(E)^\circ}\left(1-\e^{-\nu(f)}\right)
G_t^s(\d\nu)\bigg\}. \label{4.10}
 \eeqlb
Then we appeal the first equation in (\ref{4.8}) to see that
 \beqlb
\gamma_t &=&
\int_0^{t-r}\eta_s\d s + \int_0^r\d s\int_E\eta_s(\d x)
\lambda_{t-r}(x,\cdot) \nonumber \\
&=& \gamma_{t-r} + \int_E\gamma_r(\d x)\lambda_{t-r}(x,\cdot).
\label{4.11}
 \eeqlb
Combining (\ref{4.9}), (\ref{4.10}) and (\ref{4.11}) we get
 \beqnn
&&\bE\left[\exp\{- Y_t(f)\}\big| Y_s: 0\le s\le r\right] \\
&&\qquad = \exp\bigg\{-Y_r(V_{t-r}f) - \gamma_{t-r}(f) -
\int_{M(E)^\circ} \left(1-\e^{-\nu(f)}\right)
G_{t-r}(\d\nu)\bigg\},
 \eeqnn
that is, $\{Y_t: t\ge0\}$ is a Markov process with transition
semigroup $(Q^N_t)_{t\ge0}$. The theorem is proved. \qed

We next consider the semigroup $(Q_t)_{t\ge0}$ of the
$(\xi,\phi)$-superprocess. Let ${\cal K}^1(Q)$ denote the set of
probability entrance laws $K = (K_t)_{t>0}$ for the semigroup
$(Q_t)_{t\ge0}$ such that
 \beqlb
\int_0^1\d s\int_{M(E)^\circ}\nu(E)K_s(\d\nu) < \infty.
\label{4.12}
 \eeqlb
Let ${\cal K}(P)$ be the set of entrance laws $\kappa=
(\kappa_t)_{t>0}$ for the underlying semigroup $(P_t)_{t\ge0}$
that satisfy $\int_0^1\kappa_s(E)\d s < \infty$. For $\kappa\in
{\cal K}(P)$, set
 \beqlb
S_t(\kappa,f) = \kappa_t(f) - \int_0^t\d s\int_E \phi(y,V_sf(y))
\kappa_{t-s}(\d y), \quad t>0,\ f\in B(E)^+. \label{4.13}
 \eeqlb
Note that $S_t(\kappa,f) = \mu(V_tf)$ if $(\kappa_t)_{t>0}$ is
given by $\kappa_t = \mu P_t$. The following theorem characterizes
completely the set of infinitely divisible probability entrance
laws for $(Q_t)_{t\ge0}$.

\proclaim{\bf Theorem\n{4.4} {\rm (Li \cite{L96b})}.} {\sl Any
$K\in{\cal K}^1(Q)$ is infinitely divisible if and only if it is
given by
 \beqlb
\log L_{K_t}(f) = - S_t(\kappa,f) - \int_{{\cal K}(P)}
\left(1-\exp\left\{-S_t(\eta,f)\right\}\right)F(\d\eta),
\label{4.14}
 \eeqlb
where $\kappa\in{\cal K}(P)$ and $F$ is a $\sigma$-finite measure
on ${\cal K}(P)$ satisfying
 \beqlb
\int_0^1 \d s\int_{{\cal K}(P)} \eta_s(1) F(\d\eta) <\infty.
\label{4.15}
 \eeqlb
} \endproclaim
\qed

Let ${\cal K}(Q^\circ)$ be the set of entrance laws $K$ for
$(Q_t^\circ)_{t\ge0}$ satisfying (\ref{4.12}). We can also give a
general characterization for ${\cal K}(Q^\circ)$ as follows. See
also Dynkin \cite{D89}.

\proclaim{\bf Theorem\n{4.5}.} {\sl Any $H\in {\cal K}(Q^\circ)$
can be represented as
 \beqlb
& &\int_{M(E)^\circ} \left(1-\e^{-\nu(f)}\right)H_t(\d\nu) \nonumber \\
&=& S_t(\kappa,f) + \int_{{\cal K}(P)}
\left(1-\exp\left\{-S_t(\eta,f)\right\}\right) F(\d\eta), \quad
t>0,\ f\in B(E)^+, \label{4.16}
 \eeqlb
where $\kappa\in{\cal K}(P)$ and $F$ is a $\sigma$-finite measure
on ${\cal K}(P)$ satisfying (\ref{4.15}). If, in addition,
 \beqlb
\int_a^\infty\big[\sup_{x\in E}|\phi(x,z)^{-1}|\big]\d z<\infty
\label{4.17}
 \eeqlb
for some constant $a>0$, then (\ref{4.16}) defines an entrance law
$H\in {\cal K}(Q^\circ)$ for any $\kappa\in{\cal K}(P)$ and any
$\sigma$-finite measure $F$ on ${\cal K}(P)$ satisfying
(\ref{4.15}).}
\endproclaim

\noindent{\it Proof.} If $H\in {\cal K}(Q^\circ)$, then $(K)_{t>0}
= I(0,H_t)_{t>0}$ defines an infinitely divisible probability
entrance law $K\in {\cal K}^1(Q)$. Thus the representation
(\ref{4.16}) follows by (\ref{4.14}). If (\ref{4.17}) holds, there
is a family of $\sigma$-finite measures $\{L_t(x,\cdot): t>0, x\in
E\}$ on $M(E)^\circ$ such that $Q_t(\delta_x,\cdot) =
I(0,L_t(x,\cdot))$; see Dawson \cite[pp195-196] {D93}. Using this
one can show that an arbitrary infinitely divisible probability
entrance law $K\in {\cal K}^1(Q)$ may be given as $(K)_{t>0} =
I(0,H_t)_{t>0}$ for some $H\in {\cal K}(Q^\circ)$. From
(\ref{4.14}) we know that (\ref{4.16}) defines the entrance law
$H\in {\cal K}(Q^\circ)$. \qed

Let $H\in{\cal K}(Q^\circ)$ and let $\bQ^H$ be corresponding the
Kuznetsov measure supported by $W_0(M(E))$. If $N(\d s,\d w)$ is a
Poisson random measure on $[0,\infty)\times W_0(M(E))$ with
intensity $\d s\times \bQ^H(\d w)$, then
 \beqlb
Y_t = \int_{[0,t)}\int_{W_0(M(E))} w_{t-s}\, N(\d s,\d w),
\quad t\ge0,
\label{4.18}
 \eeqlb
defines an immigration process corresponding to the SC-semigroup
$(N_t)_{t\ge0}$ given by
 \beqlb
\log L_{N_t}(f) = - \int_0^t\d s\int_{M(E)^\circ}
\left(1-\e^{-\nu(f)}\right)H_s(\d\nu), \quad t\ge 0,\ f\in B(E)^+.
\label{4.19}
 \eeqlb
Clearly, (\ref{4.18}) is essentially a special form of
(\ref{4.4}). This construction has been used in \cite{L96b},
\cite{LS95} and \cite{S90}. It is simple to see from
Theorems\n{3.7},\n{4.4} and\n{4.5} that, under condition
(\ref{4.17}), any homogeneous immigration process associated with
the $(\xi,\phi)$-superprocess can be constructed in the form
(\ref{4.18}). Another related work is Evans \cite{E93}, where a
conditioned $(\xi,\phi)$-superprocess was constructed by adding up
masses thrown off by an ``immortal particle'' moving around as a
copy of $\xi$.

The construction using Kuznetsov process makes it possible to
generalize some existing results for $(\xi,\phi)$-superprocess to
the immigration process. As an example, let us give a
characterization for the ``weighted occupation time'' of the
immigration process by using the construction (\ref{4.18}). For
simplicity we only consider a special case. It is known that if
$X$ is a $(\xi,\phi)$-superprocess, then
 \beqlb
\bQ_\mu\exp\bigg\{-X_t(f)-\int_0^t X_s(g)\d s\bigg\}
=\exp\left\{-\mu(V_t(f,g))\right\}, \quad f,g\in B(E)^+,
\label{4.20}
 \eeqlb
where $V_t(f,g)(x) \equiv u_t(x)$ is the solution to
 \beqlb
u_t(x) + \int_0^t\d s\int_E \phi(x,u_s(y)) P_{t-s}(x,\d y)
=
P_tf(x) + \int_0^tP_sg(x)\d s,
\quad t\ge 0;
\label{4.21}
 \eeqlb
see e.g. Fitzsimmons \cite{F88} and Iscoe \cite{I86}. The formulas
(\ref{4.20}) and (\ref{4.21}) characterize the joint distribution
of $X_t$ and the weighted occupation time $\int_0^tX_s\d s$. By
Theorems\n{3.7} and\n{4.4} we know that
 \beqlb
\int_{M(E)} \e^{-\nu(f)} Q^\kappa_t(\mu,\d\nu) =
\exp\bigg\{-\mu(V_tf) -\int_0^t S_r(\kappa,f)\d r\bigg\}, \quad
t\ge 0, f\in B(E)^+, \label{4.22}
 \eeqlb
defines the transition semigroup $(Q^\kappa_t)_{t\ge0}$ of an
immigration process associated with the $(\xi,\phi)$-superprocess.
Let $h= \int_0^1P_s1\d s$ $\in B(E)^+$. From the discussions in
\cite{L96b} we know that $(Q_t^\kappa)_{t\ge0}$ has a realization
$(W, {\cal G}, {\cal G}_t, Y_t, \bQ_\mu^\kappa)$ such that for any
$g\in B(E)^+$ the path $\{Y_t(g\land h): t\ge0\}$ is a.s.
measurable and locally bounded, hence $\int_0^t Y_s(g)\d s$ can be
defined a.s. by increasing limits.

\proclaim{\bf Theorem\n{4.6}.} {\sl Suppose that condition
(\ref{4.17}) holds. Let $(W, {\cal G}, {\cal G}_t, Y_t,
\bQ_\mu^\kappa)$ be the realization of $(Q_t^\kappa)_{t\ge0}$
described above. Then we have
 \beqnn
&&\bQ^{\kappa}_\mu\exp\bigg\{-Y_t(f)-\int_0^t Y_s(g)\d s\bigg\} \\
&&\qquad = \exp\bigg\{-\mu(u_t)-\int_0^t S_r(\kappa,f,g) \d
r\bigg\}, \quad f,g\in B(E)^+,
 \eeqnn
where $u_t(x)$ is defined by (\ref{4.21}) and
 \beqlb
S_t(\kappa,f,g) = \kappa_t(f) + \int_0^t\kappa_s(g)\d s -
\int_0^t\kappa_{t-s}(\phi(u_s))\d s, \quad t>0. \label{4.23}
 \eeqlb
} \endproclaim

\noindent{\it Proof.} By Theorem\n{4.5} we have an entrance law
$H\in {\cal K}(Q^\circ)$ such that
 \beqlb
\int_{M(E)^\circ}\left(1-\e^{-\nu(f)}\right) H_t(\d\nu) =
S_t(\kappa,f), \quad t>0,\ f\in B(E)^+. \label{4.24}
 \eeqlb
Let $\bQ^H$ be corresponding the Kuznetsov measure on $W(M(E))$.
By (\ref{4.24}) we have
 \beqnn
& &\bQ^H\bigg(1 - \exp\bigg\{- w_t(f) - \int_0^t w_s(g)\d s\bigg\}\bigg) \\
&=& \lim_{r\downarrow0}\bQ^H\bigg(1 -
\bQ_{w_r}\exp\bigg\{-X_{t-r}(f)
- \int_0^{t-r} X_s(g)\d s\bigg\}\bigg) \\
&=&
\lim_{r\downarrow0}\bQ^H\left(1 - \exp\left\{-w_r(u_{t-r})\right\}\right) \\
&=&
\lim_{r\downarrow0}S_r(\kappa,u_{t-r}) \\
&=& S_t(\kappa,f,g),
 \eeqnn
where we have also appealed (\ref{4.21}) and (\ref{4.23}) to get
the last equality. Then using the construction (\ref{4.18}) we get
 \beqnn
& &\bQ^{\kappa}_0\exp\bigg\{-Y_t(f) - \int_0^t Y_s(g)\d s\bigg\} \\
&=& \exp\bigg\{-\int_0^t\bQ^H\bigg(1-\exp\bigg\{-w_{t-r}(f)
- \int_r^t w_{s-r}(g)\d s\bigg\}\bigg)\d r\bigg\} \\
&=& \exp\bigg\{-\int_0^t S_{t-r}(\kappa,f,g)\d r\bigg\},
 \eeqnn
and the desired result follows by the relation $\bQ_\mu^\kappa =
\bQ_\mu * \bQ_0^\kappa$. \qed

\section{Almost sure behavior of Kuznetsov processes}
\setcounter{equation}{0} In this section we study the behavior of
Kuznetsov processes near their birth times. The discussion is of
interest in providing insights into the trajectory structures of
the immigration process. Again, the lack of Feller property makes
the proof a little bit longer than expected. Let $(Q_t)_{t\ge0}$
be the transition semigroup of a $(\xi,\phi)$-superprocess.

We shall need to consider two topologies on the space $E$: the
original topology and the Ray topology of $\xi$. We write $E_r$
for the set $E$ furnished with the Ray topology of $\xi$. The
notation $M(E_r)$ is self-explanatory. Let $(P_t^b)_{t\ge0}$ be
the semigroup of bounded kernels on $E$ defined by
 \beqlb
P_t^bf(x) = \bP_xf(\xi_t)\exp\bigg\{-\int_0^t b(\xi_s)\d s\bigg\},
\quad x\in E,\ f\in B(E)^+. \label{5.1}
 \eeqlb
It is simple to check that, for any $H\in{\cal K}(Q^\circ)$,
 \beqlb
\gamma_t(f) = \int_{M(E)^\circ} \nu(f) H_t(\d\nu), \quad t> 0,\
f\in B(E)^+, \label{5.2}
 \eeqlb
defines an entrance law $\gamma= (\gamma_t)_{t>0}$ for $(P_t^b)
_{t\ge0}$.

We first consider a special $\sigma$-finite entrance law. Recall
the general formula (\ref{4.16}). Let $x\in E$ and suppose
 \beqlb
\int_{M(E)^\circ}\left(1-\e^{-\nu(f)}\right) L_t(x,\d\nu) =
V_tf(x), \quad t> 0,\ f\in B(E)^+, \label{5.3}
 \eeqlb
defines an entrance law $L(x)\in {\cal K}(Q^\circ)$. Clearly,
$(P_t^b(x,\cdot))_{t>0}$ is a minimal entrance law for
$(P_t^b)_{t\ge0}$, which may be given by (\ref{5.2}) with
$H_t(\d\nu)$ replaced by $L_t(x,\d\nu)$. From those facts it can
be deduced easily that $L(x)\in{\cal K}(Q^\circ)$ is minimal.

\proclaim{\bf Theorem\n{5.1}.} {\sl Let $\bQ^{L(x)}$ denote the
Kuznetsov measure on $W(M(E))$ determined by $L(x)\in{\cal
K}(Q^\circ)$. Then we have $w_t(E)\to0$ and $w_t(E)^{-1}
w_t\to\delta_x$ in $M(E_r)$ as $t\downarrow0$ for
$\bQ^{L(x)}$-a.a. paths $w\in W(M(E))$.} \endproclaim

\noindent{\it Proof.} The results were proved in Li and Shiga
\cite{LS95} for the case where $(P_t)_{t\ge0}$ is Feller and
$\phi(x,z) \equiv z^2/2$ by a theorem of Perkins \cite{P92} which
asserts that a conditioned $(\xi,\phi)$-superprocess is a
generalized Fleming-Viot superprocess. The calculations in
\cite{LS95} are complicated and cannot be generalized to the
present situation. We here give a proof of the theorem based on an
$h$-transform of the $(\xi,\phi)$-superprocess. The Ray cone for
the underlying process $\xi$ plays an important role in our proof.
We shall assume that $(P_t)_{t\ge0}$ is conservative. The proof
for a non-conservative underlying semigroup can be reduced to this
case as in \cite{LS95}.

Let ${\cal R}$ be a countable Ray cone for $\xi$ as constructed in
Sharpe \cite{S88} and let $\bar E$ be the corresponding Ray-Knight
compactification of $E$ with the Ray topology. Note that each
$f\in{\cal R}$ is continuous on $E_r$ and admits a unique
continuous extension $\bar f$ to $\bar E$. We regard $M(E_r)$ as a
topological subspace of $M(\bar E)$ in the usual way. Since $\bar
E$ is a compact metric space, $M(\bar E)$ is locally compact and
separable. For any fixed $u>0$,
 \beqlb
U^r_t(\mu,\d\nu)=\mu(P^b_{u-r}1)^{-1}\nu(P^b_{u-t}1)Q_{t-r}(\mu,\d\nu),
\quad 0\le r\le t\le u,
\label{5.4}
 \eeqlb
defines an inhomogeneous transition semigroup $(U^r_t)_{r\le t}$
on $M(E)^\circ$. We define the probability measure
$\bU_u^{L(x)}(\d w)$ on $W(M(E))$ by
 \beqnn
\bU_u^{L(x)}(\d w) = P^b_u 1(x)^{-1} w_u(1)\bQ^{L(x)}(\d w).
 \eeqnn
Then $\{w_t: 0<t\le u\}$ under $\bU_u^{L(x)}$ is a Markov process
with semigroup $(U^r_t)_{r\le t}$ and one-dimensional
distributions
 \beqlb
H_t(x,\d\nu):= P^b_u 1(x)^{-1}\nu(P^b_{u-t}1) L_t(x,\d\nu),
\quad 0<t\le u.
\label{5.5}
 \eeqlb
Since $L(x)\in{\cal K}(Q^\circ)$ is minimal, $(H_t(x,\cdot))
_{0<t\le u}$ is a minimal (probability) entrance law for
$(U^r_t)_{r\le t}$. Take $f\in{\cal R}$. By (\ref{5.3}) --
(\ref{5.5}) and the martingale convergence theorem we have
$\bU_u^{L(x)}$-a.s.
 \beqlb
V_t f(x) &=&
\int_{M(E)^\circ}\left(1-\e^{-\nu(f)}\right)\nu(P^b_{u-t}1)
^{-1} H_t(x,\d\nu)P^b_u 1(x)  \nonumber \\
&=&
\lim_{r\downarrow0}\int_{M(E)^\circ}\left(1-\e^{-\nu(f)}\right)
\nu(P^b_{u-t}1)^{-1} U^r_t(w_r,\d\nu)P^b_u 1(x)  \nonumber \\
&=& \lim_{r\downarrow0}w_r(P^b_{u-r}1)^{-1}
\left(1-\exp\left\{-w_r(V_{t-r}f)\right\}\right)P^b_u 1(x).
\label{5.6}
 \eeqlb
By (\ref{5.1}) and (\ref{5.6}) it follows that $\bU_u^{L(x)}$-a.s.
 \beqnn
V_tf(x) \le \liminf_{r\downarrow0}w_r(P^b_{u-r}1)^{-1}P^b_u 1(x)
\le \liminf_{r\downarrow0}\e^{2\|b\|u}w_r(1)^{-1}.
 \eeqnn
Note that $V_tf(x)$ is right continuous in $t\ge0$. Then letting
$t\downarrow0$ and $f\uparrow\infty$ in the above inequality
yields that $\bU_u^{L(x)}$-a.s. $w_t(1)\to0$ as $t\downarrow0$.
Since for each $u>0$ the measures $\bU_u^{L(x)}$ and $\bQ^{L(x)}$
are mutually absolutely continuous on $\{w\in W_0(M(E)):
w_u(1)>0\}$, we obtain the first assertion. By the same reasoning
as (\ref{5.6}) we have $\bU_t^{L(x)}$-a.s.
 \beqlb
P^b_tf(x)
&=& \int_{M(E)^\circ}\nu(f)\nu(1)^{-1}H_t(x,\d\nu)P^b_t1(x) \nonumber \\
&=& \lim_{r\downarrow0}w_r(P^b_{t-r}1)^{-1}w_r(P^b_{t-r}f)P^b_t
1(x). \label{5.7}
 \eeqlb
Clearly, $\bU_u^{L(x)}$ is absolutely continuous relative to
$\bU_t^{L(x)}$ for $u\ge t>0$. Since $f\in{\cal R}$ is an
$\alpha$-excessive function for $(P_t)_{t\ge0}$ for some
$\alpha=\alpha(f)\ge0$, from (\ref{5.1}) and (\ref{5.7}) it
follows that $\bU_u^{L(x)}$-a.s.
 \beqnn
\e^{-\|b\|t}P_tf(x) \le
\liminf_{r\downarrow0}\e^{(3\|b\|+\alpha)t}w_r(1)^{-1}w_r(f).
 \eeqnn
Take $w\in W_0(M(E))$ along which the above inequality holds for
all $f\in{\cal R}$ and all rational $t\in(0,u]$. Let $r_k=r_k(w)$
be a sequence such that $r_k\downarrow0$ and $w_{r_k}
(1)^{-1}w_{r_k} \to\hat w_0$ in $M(\bar E)$ as $k\to\infty$, where
$\hat w_0$ is a probability measure on $\bar E$. Then we have
 \beqnn
\e^{-\|b\|t}P_tf(x) \le \e^{(3\|b\|+\alpha)t}\hat w_0(\bar f).
 \eeqnn
Letting $t\downarrow0$ gives $f(x)\le\hat w_0(\bar f)$, so we have
$\hat w_0=\delta_x$. Those clearly imply $w_r(1)^{-1}w_r
\to\delta_x$ in $M(E_r)$ as $r\downarrow0$, and the second
assertion follows immediately. \qed

Now we consider an $h$-transform of the underlying semigroup
$(P_t)_{t\ge0}$. Let $h(x)= \int_0^1 P_s1(x)\d s$ for $x\in E$.
Since $h\in B(E)^{+}$ is an excessive function for
$(P_t)_{t\ge0}$, the formula
 \beqlb
T_tf(x)
=
h(x)^{-1}\int_E f(y)h(y) P_t(x,\d y),
\quad x\in E,\ f\in B(E)^+,
\label{5.8}
 \eeqlb
defines a Borel right semigroup $(T_t)_{t\ge0}$ on $E$; see e.g.
Sharpe \cite{S88}. Let $(T^\partial_t)_{t\ge0}$ be a conservative
extension of $(T_t)_{t\ge0}$ to $E^\partial := E \cup
\{\partial\}$, where $\partial$ is the cemetery point. Let
$E^{\partial}_D$ denote the entrance space of
$(T^\partial_t)_{t\ge0}$ with the Ray topology. Let $E^T_D =
E_D^{\partial} \setminus \{\partial\}$ and let $(\bar
T_t)_{t\ge0}$ be the Ray extension of $(T^\partial_t)_{t\ge0}$ to
$E^T_D$. Then $(\bar T_t)_{t\ge0}$ is also a Borel right
semigroups. Let $\kappa\in {\cal K}(P)$ be non-trivial and assume
 \beqlb
\int_{M(E)^\circ}\left(1-\e^{-\nu(f)}\right) H_t(\d\nu) =
S_t(\kappa,f), \quad t> 0,\ f\in B(E)^+, \label{5.9}
 \eeqlb
defines an entrance law $H:= L\kappa\in {\cal K}(Q^\circ)$. Let
$\bQ^{L\kappa}$ denote the corresponding Kuznetsov measure on
$W(M(E))$. Then we have

\proclaim{\bf Theorem\n{5.2}.} {\sl For $w\in W(M(E))$ define the
$M(E^T_D)$-valued path $\{h\bar w_t: t\in\IR\}$ by
 \beqlb
h\bar w_t(E^T_D\setminus E)=0 \mbox{ and }h\bar w_t(\d x)=
h(x)w_t(\d x) \mbox{ for } x\in E. \label{5.10}
 \eeqlb
Then for $\bQ^{L\kappa}$-a.a. $w\in W(M(E))$, $\{h\bar w_t: t>0\}$
is right continuous in the topology of $M(E^T_D)$ and $h\bar
w_t\to0$ as $t\downarrow0$. Moreover, for $\bQ^{L\kappa}$-a.a.
$w\in W(M(E))$ we have $w_t(h)^{-1}h\bar w_t\to\delta_{x(w)}$ for
some $x(w)\in E^T_D$ as $t\downarrow0$.} \endproclaim

\noindent{\it Proof.} By the results in Fitzsimmons \cite{F88}, if
$f\in B(E)$ is finely continuous relative to $(P_t)_{t\ge0}$, then
$\{w_t(f): t>0\}$ is right continuous for a.a. $w\in W(M(E))$.
Since the excessive function $h\in B(E)^{+}$ is finely continuous,
so is $fh$ for any bounded continuous function $f$ on $E$. It
follows that $\{hw_t: t>0\}$ is right continuous for a.a. $w\in
W(M(E))$. We may define a cumulant semigroup $(U_t)_{t\ge0}$ by
$U_tf = h^{-1}V_t(hf)$. Then $\{hw_t: t>0\}$ is a Markov process
with Borel right transition semigroup given by (\ref{1.3}) with
$(V_t)_{t\ge0}$ replaced by $(U_t)_{t\ge0}$. Let $E^T_r$ denote
the set $E$ furnished with the relative topology from $E^T_D$.
Applying the results in \cite{F88} again we conclude that $\{hw_t:
t>0\}$ is right continuous in $M(E^T_r)$ for a.a. $w\in W(M(E))$.
Therefore, $\{h\bar w_t: t>0\}$ is right continuous in $M(E^T_D)$
for a.a. $w\in W(M(E))$. Note that
 \beqnn
\bar\psi(x,z) =\left\{\begin{array}{ll}
h(x)^{-1}\phi(x,h(x)z) &\mbox{ if }x\in E, \\
0 &\mbox{ if }x\in E^T_D\setminus E,
\end{array}\right.
 \eeqnn
defines a branching mechanism $\bar\psi(\cdot,\cdot)$ on $E^T_D$.
Let $(\bar U_t)_{t\ge0}$ be the cumulant semigroup given by
 \beqlb
\bar U_t\bar f(x)
+ \int_0^t\d s\int_{E^T_D} \bar\psi(y,\bar U_s\bar f(y))
\bar T_{t-s}(x,\d y) = \bar T_t\bar f(x),
\ \ \ t\ge 0, x\in E^T_D.
\label{5.11}
 \eeqlb
Then $(\bar U_t)_{t\ge0}$ corresponds to Borel right transition
semigroup $(\bar Q_t)_{t\ge0}$ on $M(E^T_D)$. For any $t>0$ and
$x\in E^T_D$, the measure $\bar T_t(x,\cdot)$ is supported by $E$,
so $\bar T_t\bar f(x)$ and $\bar U_t\bar f(x)$ are independent of
the values of $\bar f$ on $E^T_D\setminus E$. Indeed, if $f=\bar
f|_E$ for $f\in B(E^T_D)^+$, then $\bar U_t\bar f(x) = U_tf(x)$
for all $x\in E$. We may write $\bar T_tf$ and $\bar U_tf$ instead
of $\bar U_t\bar f$ and $\bar U_t\bar f$ respectively. Clearly,
the definitions of $\bar T_tf$ and $\bar U_tf$ can be extended to
all non-negative Borel functions $f$ on $E$ by increasing limits.
As shown in \cite{L96b}, there exists a measure $\rho\in M(E^T_D)$
such that $\kappa_t(f) = \rho(\bar T_t(h^{-1}f))$ and
$S_t(\kappa,f) = \rho(\bar U_t(h^{-1}f))$. Then $\{h\bar w_t:
t>0\}$ is a Markov process with transition semigroup $(\bar
Q_t)_{t\ge0}$ and
 \beqnn
\bQ^{L\kappa} \big(1-\e^{-h\bar w_t(\bar f)}\big) = \rho(\bar
U_t\bar f), \quad t>0,\ f\in B(E^T_D)^+.
 \eeqnn
Now the results follow by Theorem\n{5.1} applied to $(\bar
U_t)_{t\ge0}$ and $(\bar T_t)_{t\ge0}$. \qed

By (\ref{4.14}), we have an entrance law $K := l\kappa\in{\cal
K}^1(Q)$ given by
 \beqlb
\int_{M(E)} \e^{-\nu(f)} K_t(\d\nu) =
\exp\left\{-S_t(\kappa,f)\right\}, \quad t>0,\ f\in B(E)^+.
\label{5.12}
 \eeqlb
It is easy to see that the restriction of $K$ to $M(E)^\circ$
belongs ${\cal K}(Q^\circ)$. Let $\bQ^{l\kappa}$ denote the
corresponding Kuznetsov measure on $W(M(E))$.

\proclaim{\bf Theorem\n{5.3}.} {\sl For $\bQ^{l\kappa}$-a.a. $w\in
W(M(E))$, $\{h\bar w_t:t>0\}$ is right continuous and $h\bar
w_t\to \rho$ for some $\rho\in M(E^T_D)$ as $t\downarrow0$.}
\endproclaim

\noindent{\it Proof.} We use the notation introduced in the proof
of Theorem\n{5.2}. Clearly, $\{h\bar w_t: t>0\}$ under
$\bQ^{l\kappa}$ is a Markov process with transition semigroup
$(\bar Q_t)_{t\ge0}$ and
 \beqnn
\bQ^{l\kappa}\exp\left\{-h\bar w_t(\bar f)\right\} =
\exp\left\{-\rho(\bar U_t\bar f)\right\}, \quad t>0,\ f\in
B(E^T_D)^+.
 \eeqnn
Thus the assertions hold by the uniqueness of transition probability. \qed

Finally, we consider the path behavior of the Kuznetsov process
determined by a general entrance rule. Let $(J_t)_{t\in\IR}$ be an
entrance rule for $(Q_t^\circ)_{t\ge0}$ satisfying
 \beqlb
\int_r^t\d s\int_{M(E)^\circ}\nu(E)J_s(\d\nu) < \infty,
\quad r\le t\in\IR.
\label{5.13}
 \eeqlb
Then we may assume that $(J_t)_{t\in\IR}$ is given by (\ref{4.2})
with the entrance laws $\{(H_{s,s+t})_{t>0}: s\in\IR\}$ taken from
${\cal K}(Q^\circ)$.

\proclaim{\bf Theorem\n{5.4}.} {\sl In the situation described
above, for $\bQ^J$-a.a. paths $w\in W(M(E))$ the process $\{h\bar
w_t: t\in\IR\}$ defined by (\ref{5.10}) is right continuous in
$M(E^T_D)^\circ$ on the interval $(\alpha(w),\beta(w))$ and $h\bar
w_t\to h\bar w_{\alpha}$ for some $h\bar w_{\alpha}\in M(E^T_D)$
as $t\downarrow\alpha(w)$. Moreover, for $\bQ^J$-a.a. paths $w\in
W(M(E))$ with $h\bar w_{\alpha}=0$, we have $w_t(h)^{-1}h \bar
w_t\to\delta_{x(w)}$ for some $x(w)\in E^T_D$ as
$t\downarrow\alpha(w)$.}
\endproclaim

\noindent{\it Proof.} Let $\bQ^H$ be the Kuznetsov measure on
$W(M(E))$ corresponding to an entrance law $H\in {\cal
K}(Q^\circ)$ represented by (\ref{4.16}). Then we have
 \beqnn
\bQ^H(\d w)=\bQ^{L\kappa}(\d w) + \int_{{\cal K}(P)}\bQ^{l\eta}(\d
w)F(\d\eta), \quad w\in W(M(E)).
 \eeqnn
By Theorems\n{5.2} and\n{5.3}, for $\bQ^H$-a.a. $w\in W(M(E))$ the
process $\{h\bar w_t: t>0\}$ is right continuous in $M(E^T_D)$ and
$h\bar w_t\to h\bar w_0$ for some $h\bar w_0\in M(E^T_D)$ as
$t\downarrow0$. Furthermore, for $\bQ^H$-a.a. $w\in W(M(E))$ with
$h\bar w_0=0$, we have $w_t(h)^{-1}h\bar w_t\to\delta_{x(w)}$ for
some $x(w)\in E^T_D$ as $t\downarrow0$. Then the desired result
holds by the representation (\ref{4.3}) of the measure $\bQ^J(\d
w)$. \qed

Clearly, (\ref{5.13}) is satisfied by the entrance rule
$(G_t)_{t\in\IR}$ in Theorem\n{4.3}. The following example shows
that the consideration of $\{h\bar w_t: t\ge 0\}$ is necessary if
one hopes to get the right limit of the path at $\alpha=
\alpha(w)$ in the usual sense.

\medskip
\noindent{\bf Example 5.1.} Suppose that $\xi$ is the minimal
Brownian motion in a bounded domain $D\subset \IR^d$ with smooth
boundary $\partial D$. We also use $\partial$ to denote the inward
normal derivative operator at $\partial D$. For any $\gamma\in
M(\partial D)$, define $(G_t)_{t>0}$ by
 \beqnn
\int_{M(D)^\circ}\left(1-\e^{-\nu(f)}\right)G_t(\d\nu) =
\int_0^t\left(1-\exp\left\{-\gamma(\partial
V_{t-s}f)\right\}\right) \d s, \quad f\in B(D)^+,
 \eeqnn
and set $G_t=0$ for $t\le 0$. Then $(G_t)_{t\in \IR}$ form an
entrance rule for the $(\xi,\phi)$-superprocess. By simple
modifications of the proofs of Theorems\n{5.3} and\n{5.4} one may
see that $w_{\alpha+} (D) =\infty$ and $w_{\alpha+}(K) =0$ for all
compact sets $K\subset D$ and $\bQ^G$-a.a. paths $w\in W(M(D))$.
\smallskip

\section{Stationary immigration processes}

\setcounter{equation}{0}

The immigration processes formulated by SC-semigroups are closely
related with the theory of excessive measures; see e.g.
Fitzsimmons and Maisonneuve \cite{FM86}, Getoor \cite{G90} and
Dellacherie et al \cite{DMM92}. In this section, we give
formulations of some results on excessive measures in terms of
stationary immigration processes.

We first consider the semigroup $(Q_t)_{t\ge0}$ of a general Borel
right MB-process. Given two probability measures $F_1$ and $F_2$
on $M(E)$, we write $F_1\preceq F_2$ if there is some probability
$G$ such that $F_1*G=F_2$. Let ${\cal E}^*(Q)$ denote the set of
all probability measures $F$ on $M(E)$ such that
 \beqlb
\int_{M(E)^\circ}\nu(1)F(\d\nu)<\infty
\label{6.1}
 \eeqlb
and $FQ_t\preceq F$ for all $t\ge0$. We write $F\in {\cal
E}_i^*(Q)$ if $F\in {\cal E}^*(Q)$ is a stationary distribution
for $(Q_t)_{t\ge0}$, and write $F\in {\cal E}_p^*(Q)$ if $F\in
{\cal E}^*(Q)$ and $\lim_{t\to\infty} FQ_t = \delta_0$. Clearly,
we have $\delta_0\in {\cal E}_i^*(Q)$, but there can be other
non-trivial stationary distributions although we are considering
the state space $M(E)$.

Let ${\cal E}(Q^\circ)$ denote the class of all excessive measures
$F$ for $(Q_t^\circ)_{t\ge0}$ satisfying (\ref{6.1}). Let ${\cal
E}_i(Q^\circ)$ be the subset of ${\cal E}(Q^\circ)$ comprising
invariant measures, and ${\cal E}_p(Q^\circ)$ the subset of purely
excessive measures. The classes ${\cal E}(Q^\circ)$ and ${\cal
E}^*(Q)$ are closely related. Indeed, $F\in{\cal E}^*(Q)$ is
infinitely divisible if and only if $F = I(\rho,J)$ for $\rho\in
M(E)$ and $J\in{\cal E}(Q^\circ)$ satisfying
 \beqlb
\int_E\rho(\d x)\lambda_t(x,\cdot) \le \rho \quad\mbox{ and }\quad
\int_E \rho(\d x)L_t(x,\cdot) + JQ_t^\circ \le J. \label{6.2}
 \eeqlb
Under the condition (\ref{4.17}), $F\in{\cal E}^*(Q)$ is
infinitely divisible if and only if $F = I(0,J)$ for some
$J\in{\cal E}(Q^\circ)$.

The following theorem shows that ${\cal E}^*(Q)$ is identical with
the totality of stationary distributions of immigration processes
associated with $X$.

\proclaim{\bf Theorem\n{6.1}.} {\sl Let $F\in {\cal E}^*(Q)$. Then
it may be written uniquely as $F = F_i*F_p$, where $F_i =
\lim_{t\to\infty}FQ_t\in {\cal E}_i^*(Q)$ and $F_p\in {\cal
E}_p^*(Q)$. Moreover, there is a unique SC-semigroup
$(N_t)_{t\ge0}$ such that $\lim_{t\to\infty}N_t = F_p$.}
\endproclaim

\noindent{\it Proof.} Let $(N_t)_{t\ge0}$ be the distributions on
$M(E)$ satisfying $F = (FQ_t) * N_t$. By the branching property of
the semigroup $(Q_t)_{t\ge0}$ one checks for any $r\ge0$ and
$t\ge0$,
 \beqlb
(FQ_{r+t})*N_{r+t}
&=&\,
F = (FQ_t)*N_t = \{[(FQ_r)*N_r]Q_t\}*N_t \nonumber \\
&=&
(FQ_{r+t})*(N_rQ_t)*N_t.
\label{6.3}
 \eeqlb
It follows that $(N_t)_{t\ge0}$ satisfies the relation
(\ref{1.4}), so it is an SC-semigroup associated with
$(Q_t)_{t\ge0}$. By the definition of ${\cal E}^*(Q)$, we have
$FQ_{r+t}\preceq FQ_t$, so the following limits exist and give the
Laplace functionals of two probability measures $F_i$ and $F_p$:
 \beqlb
L_{F_i}(f) = \uparrow\lim_{t\uparrow\infty}L_{FQ_t}(f), \quad
L_{F_p}(f) = \downarrow\lim_{t\uparrow\infty}L_{N_t}(f), \quad
f\in B(E)^+. \label{6.4}
 \eeqlb
Clearly, $F_i\in {\cal E}_i^*(Q)$ and $F = F_i * F_p$. On the
other hand,
 \beqnn
F_i*F_p = F = (FQ_t)*N_t = F_i*(F_pQ_t)*N_t,
 \eeqnn
so $F_p = (F_pQ_t)*N_t$. Therefore $F_p\in {\cal E}^*(Q)$ and
$\lim_{t\to\infty} F_pQ_t = \delta_0$. The uniqueness of the
decomposition is immediate. \qed

It is well-known that any $J\in{\cal E}(Q^\circ)$ has the Riesz
type decomposition $J = J_i + J_p$, where $J_i\in{\cal
E}_i(Q^\circ)$ and $J_p\in{\cal E}_p(Q^\circ)$ may be represented
as $J_p = \int_0^\infty H_t\,\d t$ for some $H\in {\cal
K}(Q^\circ)$. Let $\bQ^J$ be the Kuznetsov measure on $W(M(E))$
determined by $J$ and let $N^J(\d w)$ be a Poisson random measure
with intensity $\bQ^J(\d w)$. By Lemma\n{4.1},
 \beqlb
Y_t^J := \int_{W(M(E))} w_t\, N^J(\d w),
\quad t\in\IR,
\label{6.5}
 \eeqlb
defines a stationary immigration process with one-dimensional
distribution $I(0,J)$ which corresponds to the SC-semigroup
$(N_t)_{t\ge0}$ given by (\ref{4.19}). The Kuznetsov measures
determined by $J_i$ and $J_p$ are restrictions of $\bQ^J$ to
$\{w\in W(M(E)): \alpha(w)=-\infty\}$ and $\{w\in W(M(E)):
\alpha(w)>-\infty\}$, respectively; see Fitzsimmons and
Maisonneuve \cite{FM86}. It follows that
 \beqnn
Y_t^{(p)} = \int_{W(M(E))} w_t 1_{\{\alpha>-\infty\}} N^J(\d w),
\quad t\in\IR,
 \eeqnn
defines a stationary immigration process having one-dimensional
distribution $I(0,J_p)$ and SC-semigroup $(N_t)_{t\ge0}$.
Intuitively, $\{Y_t^{(p)}: t\in\IR\}$ is the ``purely
immigrative'' part of the population. On the contrary,
 \beqnn
Y_t^{(i)} = \int_{W(M(E))} w_t 1_{\{\alpha=-\infty\}} N^J(\d w),
\quad t\in\IR,
 \eeqnn
is a stationary $(\xi,\phi)$-superprocess with one-dimensional
distribution $I(0,J_i)$, which represents the ``native'' part of
the population.

The immigration process defined by (\ref{6.5}) is usually not
right continuous, but it may have right a continuous modification.
For $w\in W(M(E))$, we set
 \beqlb
w_{t+} =\left\{\begin{array}{ll}
\lim_{s\downarrow t}w_s &\mbox{if the limit exists in $M(E)$,}  \\
0 &\mbox{if the above limit does not exist in $M(E)$,}
\end{array}\right.
\label{6.6}
 \eeqlb
and define the process $\{\bar Y_t^J: t\in\IR\}$ by
 \beqlb
\bar Y_t^J = \int_{W(M(E))} w_{t+} N^J(\d w),
\quad t\in\IR,
\label{6.7}
 \eeqlb
Then $\bar Y_t^J = Y_t^J$ a.s. since $\bQ^J\{w_{t+}\neq w_t\} =
\bQ^J\{\alpha=t\} = 0$; see \cite{FM86}. In other words, $\{\bar
Y_t^J: t\in\IR\}$ is a modification for $\{Y_t^J: t\in\IR\}$.

\proclaim{\bf Theorem\n{6.2}.} {\sl Suppose $(Q_t)_{t\ge0}$ is the
transition semigroup of a $(\xi,\phi)$-superprocess. (i)~If
$J\in{\cal E}_i(Q^\circ)$, then $\{\bar Y_t^J \equiv Y_t^J:
t\in\IR\}$ is a.s. right continuous. (ii)~If $J\in {\cal
E}_p(Q^\circ)$ is a measure potential, that is,
 \beqnn
&&\int_{M(E)^\circ} \big(1-\e^{-\nu(f)}\big) J(\d\nu) \\
&&\qquad
=\int_0^\infty\d s\int_{M(E)^\circ} \big(1-\exp\{-\nu(V_sf)\}\big)G(\d\nu),
\quad f\in B(E)^+,
 \eeqnn
where $\nu(1)G(\d\nu)$ is a finite measure on $M(E)^\circ$, then
$\{\bar Y_t^J: t\in\IR\}$ is a.s. right continuous.} \endproclaim

\noindent{\it Proof.} Since (i) is simple, we only give the proof
of (ii). For $k=1,2,\cdots$ let
 \beqnn
W_k(M(E))=\{w\in W(M(E)): w_{\alpha+}(E)\ge 1/k\}.
 \eeqnn
By the results in \cite{FM86}, the path $\{w_{t+}: t\in \IR\}$ is
right continuous for $\bQ^J$-a.a. $w\in W(M(E))$ and
$\bQ^J([\cup_{k=1}^\infty W_k(M(E))]^c)=0$. Let
 \beqnn
\bar Y^{(k)}_t = \int_{W_k(M(E))} X_t(w,\cdot) 1_{\{\alpha(w)\ge
-k\}} N^J(\d w), \quad t\in\IR.
 \eeqnn
Clearly, $\{\bar Y^{(k)}_t: t\ge -k\}$ is an immigration process
corresponding to the SC-semigroup $(N^{(k)}_t)_{t\ge 0}$ given by
 \beqnn
&&\int_{M(E)}\e^{-\nu(f)} N^{(k)}_t(\d\nu) \\
&&\qquad = \exp\bigg\{ - \int_0^t \d
s\int_{M(E)}\left(1-\e^{-\nu(V_sf)}\right) 1_{\{\nu(E)\ge
1/k\}}G(\d\nu)\bigg\}.
 \eeqnn
Observe that for each $l>-k$ the process $\{\bar Y^{(k)}_t: -k\le
t\le l\}$ is a.s. a finite sum of right continuous paths and $\bar
Y^{(k)}_t\to$ $\bar Y^J_t$ increasingly as $k\to\infty$, so the
result follows as in \cite{L96b}. \qed

\proclaim{\bf Theorem\n{6.3}.} {\sl Suppose $(Q_t)_{t\ge0}$ is the
transition semigroup of a $(\xi,\phi)$-superprocess. Let $J
\in{\cal E}(Q^\circ)$ and let $\{Y_t^J: t\in\IR\}$ be defined by
(\ref{6.5}). For each $r>0$, let
 \beqnn
Y_t^{r,J} = \int_{W(M(E))} w_t 1_{\{t\ge\alpha+r\}} N(\d w), \quad
t\in\IR.
 \eeqnn
Then $\{Y_t^{r,J}: t\in\IR\}$ is an a.s. right continuous
stationary immigration process and $Y_t^{r,J} \to Y_t^J$
increasingly a.s. as $r\downarrow0$ for every $t\in\IR$.}
\endproclaim

\noindent{\it Proof.} Clearly, $JQ_r^\circ \in{\cal E}(Q^\circ)$
and
 \beqnn
JQ_r^\circ = J_i + J_pQ_r^\circ
= J_i + \int_0^\infty H_rQ_s^\circ\d s.
 \eeqnn
Using (\ref{4.1}) one may check that the Kuznetsov measure on
$W(M(E))$ determined by $JQ_r^\circ \in{\cal E}(Q^\circ)$ is the
image of $\bQ^J$ under the mapping $\{w_t: t\in\IR\}$ $\mapsto$
$\{w_t1_{\{t>\alpha+r\}}: t\in\IR\}$. It follows that
$\{Y_t^{r,J}: t\in\IR\}$ is a stationary immigration process
corresponding to $JQ_r^\circ$. By Theorem\n{6.2}, $\{Y_t^{r,J}:
t\in\IR\}$ is a.s. right continuous. The second assertion is
immediate. \qed


\noindent
\begin{thebibliography}{99}

\bibitem{AN72}
Athreya, K.B. and Ney, P.E., Branching Processes, Springer-Verlag,
New York (1972).

\bibitem{C80}
Cohn, D.L., Measure Theory, Birkh\"auser Boston, Inc., Boston, MA
(1980).

\bibitem{D92}
Dawson, D.A., Infinitely Divisible Random Measures and
Superprocesses, In: Proceedings of 1990 Workshop on Stochastic
Analysis and Related Topics in Silivri, Progress in Probability
{\bf 31} (1992) 1-129, H. K\"orezlio\u glu and A.S. \"Ust\"unel
eds., Birkh\"auser Boston, Inc., Boston, MA.

\bibitem{D93}
Dawson, D.A., Measure-valued Markov Processes, In: Ecole d'Et\'e
de Probabilit\'es de Saint-Flour XXI-1991, Lecture Notes Math.
{\bf 1541} (1993), 1-260, P.L. Hennequin ed., Springer-Verlag, New
York.

\bibitem{DI78}
Dawson, D.A. and Ivanoff, D., Branching diffusions and random
measures, In: Branching Processes, Advances in Probability and
Related Topics {\bf 5} (1978), 61-103, A. Joffe and P. Ney eds.,
Marcel Dekker, Inc., New York.

\bibitem{DMM92}
Dellacherie, C., Maisonneuve, B. and Meyer, P.A., Probabilit\'es
et Potential, Ch. 17-24, Hermann (1992).

\bibitem{D89}
Dynkin, E.B., Three classes of infinite dimensional diffusion
processes, J. Funct. Anal. {\bf 86} (1989), 75-110.

\bibitem{D91}
Dynkin, E.B., Branching particle systems and superprocesses, Ann.
Probab. {\bf 19} (1991), 1157-1194.

\bibitem{E92}
Evans, S., The entrance space of a Markov branching process
conditioned on non-extinction, Can. Math. Bull. {\bf 35} (1992),
70-74.

\bibitem{E93}
Evans, S., Two representations of conditioned superprocess,
Proceedings of Royal Society of Edinburgh {\bf 123A} (1993),
959-971.

\bibitem{F88}
Fitzsimmons, P.J., Construction and regularity of measure-valued
Markov branching processes, Israel J. Math. {\bf 64} (1988),
337-361.

\bibitem{F92}
Fitzsimmons, P.J., On the martingale problem for measure-valued
Markov branching processes, In: Seminar on Stochastic Processes
{\bf 1991} (1992), 39-51, E. Cinlar et al eds., Birkh\"auser
Boston, Inc., Boston, MA.

\bibitem{FM86}
Fitzsimmons, P.J. and Maisonneuve, B., Excessive measures and
Markov processes with random birth and death, Probab. Theory
Related Fields {\bf 18} (1986), 571-575.

\bibitem{G90}
Getoor, R.K., Excessive Measures, Birkh\"auser Boston, Inc.,
Boston, MA (1990).

\bibitem{GG87}
Getoor, R.K. and Glover, J., Constructing Markov processes with
random times of birth and death, In: Seminar on Stochastic
Processes {\bf 1986} (1987), 35-69, E. Cinlar et al eds.,
Birkh\"auser Boston, Inc., Boston, MA.

\bibitem{GL90}
Gorostiza, L.G. and Lopez-Mimbela, J.A., The multitype measure
branching process, Adv. Appl. Probab. {\bf 22}, 49-67 (1990).

\bibitem{GR90}
Gorostiza, L.G. and Roelly, S., Some properties of the multitype
measure branching process, Stochastic Process. Appl. {\bf 37}
(1990), 259-274.

\bibitem{HL99}
Hong, W.M. and Li, Z.H., A central limit theorem for super
Brownian motion with super Brownian immigration, J. Appl. Probab.
{\bf 36} (1999), 1218-1224.

\bibitem{I86}
Iscoe, I., A weighted occupation time for a class of
measure-valued branching processes, Probab. Theory Related Fields
{\bf 71} (1986), 85-116.

\bibitem{K75}
Kallenberg, O., Random measures, Academic Press, New York (1975).

\bibitem{KW71}
Kawazu, K. and Watanabe, S., Branching processes with immigration
and related limit theorems, Theory Probab. Appl. {\bf 16} (1971),
36-54.

\bibitem{K74}
Kuznetsov, S.E., Construction of Markov processes with random
times of birth and death, Theory Probab. Appl. {\bf 18} (1974),
571-575.

\bibitem{L92a}
Li, Z.H., A note on the multitype measure branching process, Adv.
Appl. Probab. {\bf 24} (1992), 496-498.

\bibitem{L92b}
Li, Z.H., Measure-valued branching processes with immigration,
Stochastic Process. Appl. {\bf 43} (1992), 249-264.

\bibitem{L96a}
Li, Z.H., Convolution semigroups associated with measure-valued
branching processes, Chinese Sci. Bull. (Chinese Edition) {\bf
40}, 2018-2021 / (English Edition) {\bf 41} (1996), 276-280.

\bibitem{L96b}
Li, Z.H., Immigration structures associated with Dawson-Watanabe
superprocesses, Stochastic Process. Appl. {\bf 62} (1996), 73-86.

\bibitem{L98a}
Li, Z.H., Immigration processes associated with branching
particle systems, Adv. Appl. Probab. {\bf 30} (1998), 657-675.

\bibitem{L98b}
Li, Z.H., Entrance laws for Dawson-Watanabe superprocesses with
non-local branching, Acta Mathematica Scientia (Series A, English
Edition) {\bf 18} (1998), 449-456.

\bibitem{LS95}
Li, Z.H. and Shiga, T., Measure-valued branching diffusions:
immigrations, excursions and limit theorems, J. Math. Kyoto Univ.
{\bf 35} (1995), 233-274.

\bibitem{P92}
Perkins, E.A., Conditional Dawson-Watanabe processes and
Fleming-Viot processes, In: Seminar on Stochastic Processes {\bf
1991} (1992), 143-156, E. Cinlar et al eds., Birkh\"auser Boston,
Inc., Boston, MA.

\bibitem{S88}
Sharpe, M.J., General Theory of Markov Processes, Academic Press,
New York (1988).

\bibitem{S90}
Shiga, T., A stochastic equation based on a Poisson system for a
class of measure-valued diffusion processes, J. Math. Kyoto Univ.
{\bf 30} (1990), 245-279.

\bibitem{SW73}
Shiga, T. and Watanabe, S., Bessel diffusions as a one-parameter
family of diffusion processes, Z. Wahrsch. verw. Geb. {\bf 27}
(1973), 37-46.

\bibitem{S69}
Silverstein, M.L., Continuous state branching semigroups, Z.
Wahrsch. verw. Geb. {\bf 9} (1969), 235-257.

\bibitem{W68}
Watanabe, S., A limit theorem of branching processes and
continuous state branching processes, J. Math. Kyoto Univ. {\bf 8}
(1968), 141-167.

\end{thebibliography}
\end{document}